\documentclass[12pt]{article}
\usepackage{amsfonts}

\usepackage{latexsym}
\usepackage{amssymb}
\usepackage{epsfig}
\usepackage{amsmath}
\usepackage{amsthm}
\usepackage[mathscr]{eucal}
\usepackage{subfigure}
\usepackage{graphicx}

\newtheorem{Lemma}{Lemma}

\newtheorem{Theorem}[Lemma]{Theorem}

\newtheorem{Definition}{Definition}

\renewcommand{\qed}{\hfill{\ \ \rule{2mm}{2mm}} \vspace{0.2in}}

\begin{document}

\title{Duality in percolation via outermost boundaries I: Bond Percolation}
\author{ \textbf{Ghurumuruhan Ganesan}
\thanks{E-Mail: \texttt{gganesan82@gmail.com} } \\
\ \\
New York University, Abu Dhabi}
\date{}
\maketitle

\begin{abstract}



Tile \(\mathbb{R}^2\) into disjoint unit squares \(\{S_k\}_{k \geq 0}\) with the origin being the centre of \(S_0\)
and say that \(S_i\) and \(S_j\) are star adjacent if they share a corner and plus adjacent if they share an edge.
Every square is either vacant or occupied. Outermost boundaries of finite star and plus connected components
frequently arise in the context of contour analysis in percolation and random graphs.
In this paper, we derive the outermost boundaries for finite star and plus connected components using
a piecewise cycle merging algorithm.
For plus connected components, the outermost boundary is a single cycle and
for star connected components, we obtain that the outermost boundary is a connected union
of cycles with mutually disjoint interiors. As an application, we use
the outermost boundaries to give an alternate proof of mutual exclusivity of
left right and top bottom crossings in oriented and unoriented bond percolation.

\vspace{0.1in} \noindent \textbf{Key words:} Star and plus connected components, outermost boundary, union of cycles, left right and top bottom crossings.

\vspace{0.1in} \noindent \textbf{AMS 2000 Subject Classification:} Primary:
60J10, 60K35; Secondary: 60C05, 62E10, 90B15, 91D30.
\end{abstract}

\bigskip

\renewcommand{\theequation}{\thesection.\arabic{equation}}
\setcounter{equation}{0}
\section{Introduction} \label{intro}

Tile \(\mathbb{R}^2\) into disjoint unit squares \(\{S_k\}_{k \geq 0}\) with origin being the centre of~\(S_0.\) Every square in \(\{S_k\}\) is assigned one of the two states, occupied or vacant and the square \(S_0\) containing the origin is always occupied. For \(i \neq j,\) we say that~\(S_i\) and~\(S_j\) are \emph{adjacent} or \emph{star adjacent} if they share a corner between them. We say that \(S_i\) and \(S_j\) are \emph{plus adjacent}, if they share an edge between them. Here we follow the notation of Penrose (2003).

The structure of the outermost boundary of finite components is crucial for contour analysis problems of percolation (Grimmett~(1999), Ganesan~(2014)) and random graphs (Penrose~(2003), Ganesan~(2013)). For plus connected components, the boundaries have been well studied before (see for example, Penrose (2003)) and it is also relatively easy to visualize that the boundary must be a single cycle. For star connected components, visualization is a bit difficult since there are many possible candidates for the outermost boundary. Below, we give a definition of outermost boundary that holds for both star and plus connected components and derive the structure of the outermost boundary for the star connected component. The corresponding result for the plus connected component is then obtained as a corollary.

Before we present our results, we briefly enumerate some recent literature containing the related duality problem of percolation. Tim\'ar~(2013) uses separating sets in equivalence class of infinite paths to study duality in slightly more general locally finite graphs and Penrose (2003) uses unicoherence and topological arguments to investigate plus connected components. Bollob\'as and Riordan (2006) use a step by step construction for obtaining the outermost boundary in a slightly related problem in bond percolation. In essence, most of the proofs use either substantial topology or infinite graphs.

Our aim in this paper is two fold. We derive the outermost boundary structure for \emph{both} star and plus connected components. We then use the finite graph theoretic structure of the outermost boundaries to provide an alternate proof of mutual exclusivity between left right and top bottom crossings in both unoriented and oriented crossings in percolation.

\subsection*{Model Description}
We first discuss star connected components. We say that the square~\(S_i\) is connected to the square~\(S_j\) by a \emph{star connected \(S-\)path} if there is a sequence of distinct squares \((Y_1,Y_2,...,Y_t), Y_l \subset \{S_k\}, 1 \leq l \leq t\) such that~\(Y_l\) is star adjacent to~\(Y_{l+1}\) for all \(1 \leq l \leq t-1\) and \(Y_1 = S_i\) and \(Y_t = S_j.\) If all the squares in \(\{Y_l\}_{1 \leq l \leq t}\) are occupied, we say that~\(S_i\) is connected to~\(S_j\) by an \emph{occupied} star connected \(S-\)path.

Let \(C(0)\) be the collection of all occupied squares in \(\{S_k\}\) each of which is connected to the square~\(S_0\) by an occupied star connected \(S-\)path. We say that \(C(0)\) is the  star connected occupied component containing the origin. Throughout we assume that \(C(0)\) is finite and let \(\{J_k\}_{1 \leq k \leq M} \subset \{S_j\}\) be the set of all the occupied squares  belonging to the component \(C(0).\) In this paper, we study the outermost boundary for finite star connected components containing the origin and by translation, the results hold for arbitrary finite star connected components.



It is possible that a star connected component has multiple choices for the outermost boundary. Consider for example, the component consisting of the union of four squares \(S_a, S_b, S_c\) and~\(S_d\) each sharing exactly one edge with the square \(S_0\) containing the origin. If we consider each square itself as a cycle, then the union of the edges of these four cycles could itself be considered as the boundary. On the other hand, it is also possible to consider the ``bigger" cycle consisting of edges in \((S_a \cup S_b \cup S_c \cup S_d) \setminus S_0\) as the boundary.

To avoid ambiguities as described in the above paragraph, we give below a formal definition of outermost boundary. We first have a few preliminary definitions. Let \(G_0\) be the graph with vertex set being the set of all corners of the squares of \(\{S_k\}\) in the component~\(C(0)\) and edge set consisting of the edges of the squares of \(\{S_k\}\) in~\(C(0).\) Two vertices in the graph~\(G_0\) are said to be adjacent if they share an edge between them. Two edges in \(G_0\) are said to be adjacent if they share an endvertex between them.

Let \(P = (e_1,e_2,\ldots,e_t)\) be a sequence of distinct edges in \(G_0.\) We say that~\(P\) is a \emph{path} if~\(e_i\) and~\(e_{i+1}\) are adjacent for every \(1 \leq i \leq t-1.\)  Let~\(a\) be the endvertex of~\(e_1\) not common to~\(e_2\) and let~\(b\) be the endvertex of~\(e_t\) not common to~\(e_{t-1}.\) The vertices \(a\) and \(b\) are the \emph{endvertices} of the path~\(P.\)

We say that~\(P\) is a \emph{self avoiding path} if the following three statements hold: The edge \(e_1\) is adjacent only to \(e_2\) and no other~\(e_j, j \neq 2.\) The edge~\(e_t\) is adjacent only to~\(e_{t-1}\) and no other~\(e_j, j \neq t-1.\) For each~\(1 \leq i \leq t-1,\) the edge~\(e_i\) shares one endvertex with~\(e_{i-1}\) and another endvertex with~\(e_{i+1}\) and is not adjacent to any other edge~\(e_j, j\neq i-1,i+1.\)

\begin{figure}[tbp]
\centering
\includegraphics[width=2.5in, trim= 100 350 200 200, clip=true]{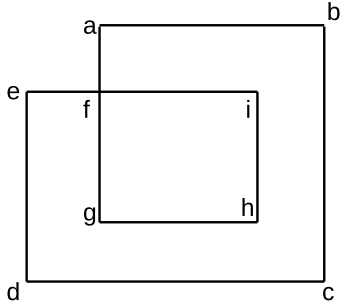}
\caption{The sequence of vertices \(abcdefghifa\) form a circuit but not a cycle.}
\label{circ_not_cyc_fig}
\end{figure}

We say that~\(P\) is a \emph{circuit} if \((e_1,e_2,\ldots,e_{t-1})\) forms a path and the edge~\(e_t\) shares one endvertex with \(e_1\) and another endvertex with~\(e_{t-1}.\) We say that~\(P\) is a \emph{cycle} if \((e_1,e_2,...,e_{t-1})\) is a self avoiding path and the edge~\(e_t\) shares one endvertex with~\(e_1\) and another endvertex with~\(e_{t-1}\) and does not share an endvertex with any other edge~\(e_j, 2 \leq j \leq t-2.\) We emphasize here that we consider only cycles that do not intersect themselves. In other words, every vertex in a cycle~\(C\) is adjacent to \emph{exactly} two edges of~\(C.\) For example, the sequence of edges formed by the vertices~\(abcdefghifa\) in Figure~\ref{circ_not_cyc_fig} is a circuit but not a cycle. The sequence of vertices \(abcdefa\) forms a cycle.

Any cycle~\(C\) contains at least four edges and divides the plane~\(\mathbb{R}^2\) into two disjoint connected regions. As in Bollob\'as and Riordan~(2006), we denote the bounded region to be the \emph{interior} of~\(C\) and the unbounded region to be the \emph{exterior} of~\(C.\) We use cycles to define the outermost boundary of star connected components.

Let \(e\) be an edge in the graph~\(G_0\) defined above. We say that~\(e\) is adjacent to a square~\(S_k\) if it is one of the edges of~\(S_k.\) We say that the edge~\(e\) is contained in the \emph{interior} of the cycle~\(C\) if both the squares in~\(\{S_k\}\) containing~\(e\) as an edge, lie in the interior of~\(C.\) An analogous definition holds for edges in the exterior of~\(C.\) We say that~\(e\) is a  \emph{boundary edge} if it is adjacent to a vacant square and is also adjacent to an occupied square of the component~\(C(0).\) Let \(C\) be any cycle of edges in \(G_0.\)
We have the following definition.
\begin{Definition} \label{out_def} We say that the edge \(e\) in the graph~\(G_0\) is an \emph{outermost boundary} edge of the component \(C(0)\) if the following holds true for every cycle \(C\) in~\(G_0:\) either \(e\) is an edge in \(C\) or \(e\) belongs to the exterior of \(C.\)

We define the outermost boundary \(\partial _0\) of \(C(0)\) to be the set of all outermost boundary edges of~\(G_0.\)
\end{Definition}
Thus outermost boundary edges cannot be contained in the interior of any cycle in the graph~\(G_0.\) Our main result is the following.
\begin{Theorem}\label{thm3} Suppose \(C(0)\) is finite.
The outermost boundary \(\partial_0\) of \(C(0)\) is the union of a unique set of cycles \(C_1,C_2,\ldots,C_n\) in \(G_0\) with the following properties:\\
\((i)\) Every edge in \(\cup_{1 \leq i \leq n} C_i\) is an outermost boundary edge.\\
\((ii)\) The graph \(\cup_{1 \leq i \leq n}C_i\) is a connected subgraph of \(G_0.\)\\
\((iii)\) If \(i \neq j,\) the cycles \(C_i\) and \(C_j\) have disjoint interiors and have at most one vertex in common.\\
\((iv)\) Every occupied square \(J_k \in C(0)\) is contained in the interior of some cycle \(C_{j}.\)\\
\((v)\) If \(e \in C_{j}\) for some \(j,\) then \(e\) is a boundary edge adjacent to an occupied square of~\(C(0)\) contained in the
interior of~\(C_j\) and also adjacent to a vacant square lying in the exterior of all the cycles in~\(\partial_0.\)\\
Moreover, there exists a circuit \(C_{out}\) containing every edge of \(\cup_{1 \leq i \leq n} C_i.\)
\end{Theorem}
The outermost boundary \(\partial_0\) is therefore also an Eulerian graph with \(C_{out}\) denoting the corresponding Eulerian circuit (for definitions, we refer to Chapter 1, Bollob\'as (2001)). We remark that the above result also provides a more detailed justification of the statement made about the outermost boundary and the corresponding circuit in the proof of Lemma~3 of Ganesan (2013).



The proof technique of the above result also allows us to obtain the outermost boundary for plus connected components. We recall that squares \(S_i\) and \(S_j\) are \emph{plus adjacent} if they share an edge between them. We say that the square~\(S_i\) is connected to the square~\(S_j\) by a \emph{plus connected \(S-\)path} if there is a sequence of distinct squares \((Y_1,Y_2,...,Y_t), Y_l \subset \{S_k\}, 1 \leq l \leq t\) such that~\(Y_l\) is plus adjacent to \(Y_{l+1}\) for all \(1 \leq l \leq t-1\) and \(Y_1 = S_i\) and \(Y_t = S_j.\) If all the squares in \(\{Y_l\}_{1 \leq l \leq t}\) are occupied. we say that \(S_i\) is connected to \(S_j\) by an \emph{occupied} plus connected \(S-\)path.

Let \(C^+(0)\) be the collection of all occupied squares in \(\{S_k\}\) each of which is connected to the square~\(S_0\) by an occupied plus connected \(S-\)path. We say that \(C^+(0)\) is the plus connected occupied component containing the origin. Throughout we assume that \(C^+(0)\) is finite. Let \(G^+_0\) be the graph with vertex set being the set of all corners of the squares of \(\{S_k\}\) in \(C^+(0)\) and edge set consisting of the edges of the squares of \(\{S_k\}\) in~\(C^+(0).\)

Every plus connected component is also a star connected component and so the definition of outermost boundary edge in Definition~\ref{out_def} holds for the component~\(C^+(0)\) with~\(G_0\) replaced by \(G^+_0.\) We have the following result.
\begin{Theorem}\label{thm2} Suppose \(C^+(0)\) is finite. The outermost boundary \(\partial^+_0\) of \(C^+(0)\) is unique cycle in \(G^+_0\) with the following properties:\\
\((i)\) All squares of \(C^+(0)\) are contained in the interior of \(\partial^+_0.\)\\
\((ii)\) Every edge in \(\partial^+_0\) is a boundary edge adjacent to an occupied square of~\(C^+(0)\) contained in the interior of \(\partial^+_0\) and a vacant square in the exterior.
\end{Theorem}
This is in contrast to star connected components which may contain multiple cycles in the outermost boundary.

To prove Theorems~\ref{thm3} and~\ref{thm1}, we use the following intuitive result about merging cycles. Let \(G\) be the graph with vertex set being the corners of the squares \(\{S_k\}_{k \geq 0}\) and edge set being the edges of the squares \(\{S_k\}_{k \geq 0}.\)
\begin{Theorem}\label{thm1} Let \(C\) and \(D\) be cycles in the graph~\(G\) that have more than one vertex in common. There exists a unique cycle~\(E\) consisting only of edges of \(C\) and \(D\) with the following properties:\\
\((i)\) The interior of \(E\) contains the interior of both \(C\) and \(D.\)\\
\((ii)\) If an edge \(e\) belongs to \(C\) or \(D,\) then either \(e\) belongs to \(E\) or is contained in its interior.

Moreover, if \(D\) contains at least one edge in the exterior of \(C,\) then the cycle \(E\) also contains an edge of \(D\) that lies in the exterior of~\(C.\)
\end{Theorem}
The above result essentially says that if two cycles intersect at more than one point, there is an innermost cycle containing both of them in its interior. We provide an iterative piecewise algorithmic construction for obtaining the cycle \(E,\) analogous to Kesten~(1980) for crossings, in Section~\ref{pf1}.

\subsection*{Bond Percolation}
In this section, we use the structure of the outermost boundaries derived in the previous subsection to give alternate proof for mutual exclusivity of left right and top bottom crossings in oriented and unoriented bond percolation. We first discuss for unoriented bond percolation and then consider oriented bond percolation. We also remark that most of the existing approaches for obtaining the mutual exclusivity mainly use some version of interface graphs (Bollob\'as and  Riordan (2006)) that is usually obtained via a step by step procedure. Our method uses the outermost boundaries to directly obtain the presence of the closed top bottom dual crossing in the absence of open left right crossings. For more material on left right and top bottom crossings we refer to Bollob\'as and Riordan (2006).

We recall that \(G\) is the graph with vertex set as the set of corners of the squares in \(\{S_k\}.\) The edge set is the set of edges of the squares in \(\{S_k\}.\) Let~\(G_d\) be the graph obtained by shifting the graph~\(G\) by \(\left(\frac{1}{2},\frac{1}{2}\right).\) The graph~\(G_d\) tiles~\(\mathbb{R}^2\) into disjoint unit squares~\(\{W_k\}\) such that \(W_k = S_k + \left(\frac{1}{2},\frac{1}{2}\right).\) We say that the graph \(G\) is the \emph{dual graph} of \(G_d.\)

Consider bond percolation in the graph \(G_d\) where every edge is either open or closed. By construction every edge \(f\) of the graph~\(G_d\) intersects perpendicularly a unique dual edge \(e = e(f)\) of \(G.\) We say that \(e\) is open if and only if \(f\) is open. Let \(R = [0,m] \times [0,n]\) be the \(m \times n\) rectangle in the graph \(G_d\) containing exactly \(mn\) edges of \(G_d.\) Let \(e\) be an edge with endvertices \(u_1\) and \(u_2.\) We say that edge \(e\) lies in the \emph{interior} of \(R\) if the following property holds. For \(i = 1,2,\) the endvertex \(u_i\) either belongs to the boundary of \(R\) or lies in the interior of \(R.\)

Let \(R_{left}\ = \{0\} \times [0,n]\) and \(R_{right} = \{m\} \times [0,n]\) be the left and right edges of \(R,\) respectively. We say that a self avoiding path \(P = (e_1,\ldots,e_k)\) is a left right crossing for \(R\) if the following three conditions hold:\\
\((a)\) The edge~\(e_1\) contains exactly one endvertex in~\(R_{left}\) and no other edge in~\(P\) intersects~\(R_{left}.\)\\
\((b)\) The edge~\(e_k\) contains exactly one endvertex in~\(R_{right}\) and no other edge in~\(P\) intersects~\(R_{right}.\)\\
\((c)\) Every edge \(e_j, 2 \leq j \leq k-1\) have both their endvertices in~\(R.\)\\
If every edge in \(P\) is open we say that \(P\) is an open left right crossing. An analogous definition holds for top bottom and open top bottom crossings of~\(R.\)

We have a similar definition for the dual crossing. We say that the dual edge \(f \in G\) intersects the left (right) edge of \(R\) if~\(f\) intersects some edge of \(G_d\) contained in the left (right) edge of~\(R.\) We define a dual left right crossing of the rectangle \(R\) as follows. We say that a self avoiding path \(P_d = (f_1,\ldots,f_m)\) in the graph \(G\) is a \emph{dual left right crossing} for \(R\) if the edge \(f_1\) intersects \(R_{left},\) the edge \(f_m\) intersects \(R_{right}\) and every other edge \(f_j, 2 \leq j \leq m-1\) have both their endvertices in the interior of the rectangle~\(R.\)
As before, we say that \(P_d\) is an \emph{open} dual left right crossing if every edge in \(P_d\) is open. An analogous definition holds for dual top bottom and open dual top bottom crossings of~\(R.\)


We have the following result regarding left right and top bottom crossings.
\begin{Theorem}\label{thm6} One of the following two events always occurs but not both:\\
\((i)\) The rectangle \(R\) contains an open left right crossing.\\
\((ii)\) The rectangle \(R\) contains a closed dual top bottom crossing.\\
By rotating rectangles, we also have that one of the following two events always occurs but not both:\\
\((iii)\) The rectangle \(R\) contains an open left right dual crossing.\\
\((iv)\) The rectangle \(R\) contains a closed top bottom crossing.
\end{Theorem}


\subsection*{Oriented bond percolation}
\begin{figure*}
\centerline{\subfigure[]{\includegraphics[width=2.5in, trim= 0 450 300 100, clip=true]{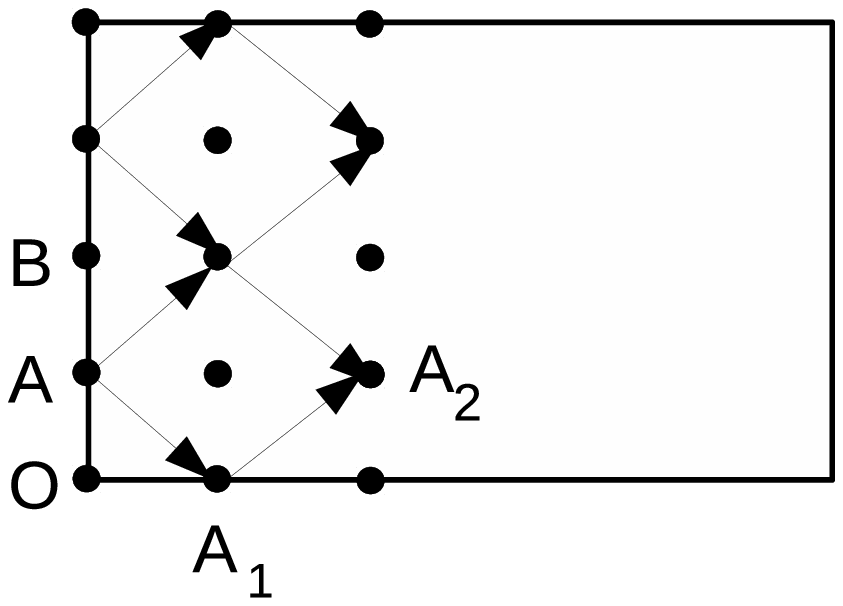}
\label{fig_first_case}}
\hfil
\subfigure[]{\includegraphics[width=2.5in, trim= 50 400 200 100, , clip=true]{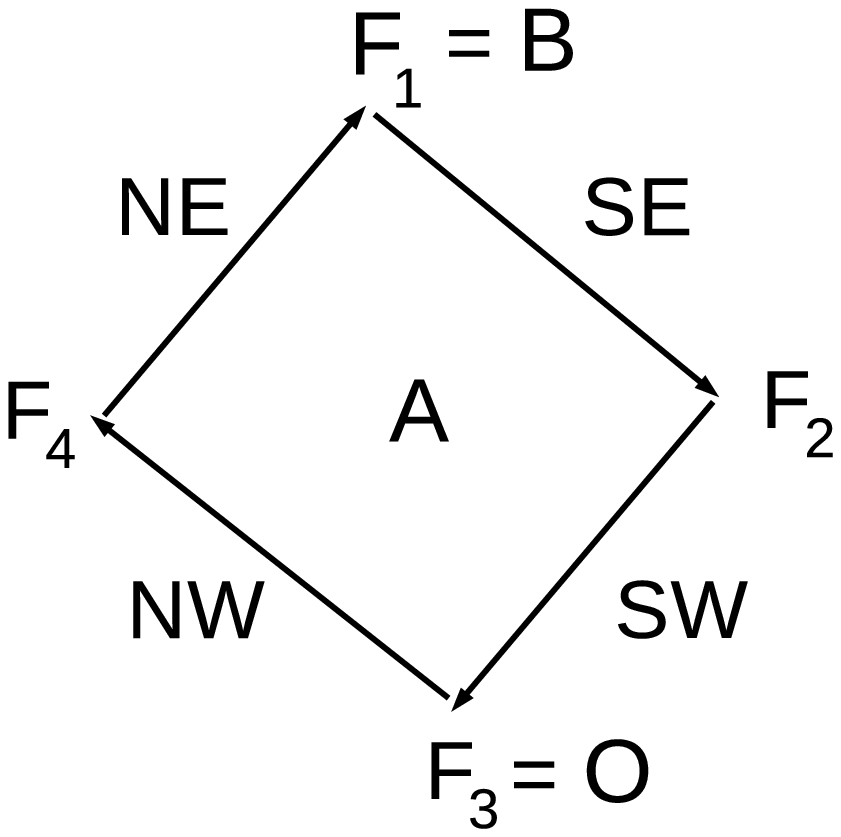}
\label{fig_second_case}}}
\caption{\((a)\) The oriented bond percolation model in the rectangle \(R.\) Here \(O\) is the origin, \(A = (0,1), B = (0,2), A_1 = (1,0)\) and \(A_2 = (2,1).\) \((b)\)~The square \(S''_A\) with centre \(A\) forming an oriented cycle is shown with the labels for the corresponding orientations.}
\label{or_mod}
\end{figure*}

As before, we consider the \(m \times n\) rectangle \(R = [0,m] \times [0,n].\) We consider the oriented edge model in \(R\) where oriented edges are present as follows.\\
\((b1)\) If \(i = 0\) or \(i\) is even, then there is an oriented edge from \((i,j)\) to \((i+1,j+1)\) and from \((i,j)\) to \((i+1,j-1),\) for \(1 \leq j \leq n-1, j\) odd. If \(n\) is odd, there is only the oriented edge from \((i,n)\) to \((i+1,n-1).\) \\
\((b2)\) If \(i\) is odd, then for \(2 \leq j \leq n-1, j\) even, there is an oriented edge from~\((i,j)\) to~\((i+1,j+1)\) and from~\((i,j)\) to~\((i+1,j-1).\) For \(j = 0,\) there is only the oriented edge from~\((i,0)\) to~\((i+1,1)\) and if \(n\) is even, there is only the oriented edge from~\((i,n)\) to~\((i+1,n-1).\)\\
We refer to Figure~\ref{or_mod}\((a)\) for illustration where we have drawn part of the oriented bond percolation model for~\(i =0\) and~\(i = 1.\)

Let \({\cal F}_{or}\) denote the set of all oriented edges with both endvertices in the rectangle \(R\) as described above. Let \(e_i \in {\cal F}_{or}, i = 1,2\) be two oriented edges in the rectangle~\(R\) and let \(x_i\) and \(y_i\) denote the non arrow and arrow endvertices of \(e_i,\) respectively. We say that \((e_1,e_2)\) is a \emph{consistent pair} if the edges share an endvertex and the arrow end \(y_1\) of~\(e_1\) coincides with the non arrow endvertex~\(x_2\) of~\(e_2;\)  i.e., \(y_1 = x_2.\)

We say that a sequence of distinct oriented edges \(P = (e_1,e_2,\ldots,e_k)\) form an \emph{oriented path} if for every \(1 \leq i \leq k-1,\) the pair of edges \((e_i,e_{i+1})\) is a consistent pair. Let~\(x\) be the endvertex of \(e_1\) not common with \(e_2\) and let~\(y\) be the endvertex of~\(e_k\) not common with~\(e_{k-1}.\) We say that \(x\) and \(y\) are the \emph{endvertices} of~\(P.\) If \(P\) is any oriented path, one endvertex of \(P\) is a non arrow endvertex and the other endvertex is an arrow endvertex. In Figure~\ref{or_mod}\((a),\) the pair of edges between the points \(AA_1\) and \(A_1A_2\) are consistent and form a oriented path with endvertices \(A\) and \(A_2.\)

Let \(R_{left} = \{0\}\times [0,n]\) and \(R_{right} = \{m\} \times [0,n]\) respectively denote the left and right edges of the rectangle \(R.\) An oriented path having one non arrow endvertex in~\(R_{left}\) and another arrow endvertex in~\(R_{right}\) is called an \emph{oriented left right crossing of \(R.\)} To every oriented edge \(e \in {\cal F}_{or},\) we assign one of the two following states: open or closed.  If every edge in an oriented left right crossing~\(P\) of~\(R\) is open, we say that~\(P\) is an open oriented left right crossing of~\(R.\)

The next step is to define the \emph{dual lattice} and we follow the notation of Durrett (1984). We tile \(\mathbb{R}^2\) into disjoint~\(1 \times 1\) squares \(\{S''_{z}\}_{z \in \mathbb{Z}^2},\) so that if \(z = (i,j),\) then \(S''_z\) has endvertices \((i,j-1),\) \((i+1,j),\) \((i,j+1)\) and \((i-1,j).\) We orient the edges of \(S''_z\) in such a way that they form a clockwise oriented cycle; i.e. an oriented path with coincident endvertices. Every edge in \(S''_z\) is therefore one of the four following types: type~\(1\)(\(\nearrow\)), type~\(2\)(\( \nwarrow\)), type~\(3\)(\(\searrow\)) and type~\(4\)(\(\swarrow\)). According to the orientations, we also call them as \(NE, NW, SE\) and~\(SW\) arrows, respectively, representing north east, north west, south east and south west directions. Edges belonging to the square~\(\{S''_z\}_{z \in R}\) are called \emph{dual} edges. In Figure~\ref{or_mod}\((b)\) we have illustrated the square \(S''_A\) with centre \(A = (0,1).\)


Dual edges can lie in the interior or the exterior of the rectangle \(R.\) Every dual edge lying in the interior of \(R\) is assigned one of the two states open or closed as follows. Let \(f\) be a dual edge contained in the interior of \(R;\) i.e., no endvertex of \(f\) lies in the exterior of \(R.\) The edge~\(f\) intersects a unique edge \(e = e(f)\) belonging to the original percolation model descsribed in Figure~\ref{or_mod}. We say that \(f\) is open if~\(e\) is open and~\(f\) is closed if~\(e\) is closed. Thus open or closed dual edges necessarily lie in the interior of \(R.\)

We say that a dual oriented path \(\Pi = (e_1,e_2,\ldots,e_k)\) is a \emph{dual oriented top bottom crossing} of the rectangle~\(R\) if the following properties \((a1)-(a3)\) hold.\\
\((a1)\) The first edge~\(e_1\) has its non-arrow endvertex in the top edge~\(R_{top}\) of~\(R.\)\\
\((a2)\) The last edge~\(e_k\) has its arrow endvertex at the bottom edge~\(R_{bottom}\) of~\(R.\)\\
\((a3)\) If \(u\) is an endvertex of an edge \(e \in \Pi,\) then either \(u\) lies on the boundary of \(R\) or lies in the interior of \(R.\)\\
If every dual edge in the path~\(\Pi\) is closed, we say that~\(\Pi\) is a \emph{closed} dual oriented top bottom crossing.


We have the following result.
\begin{Theorem}\label{thm8} One of the following two events always occurs but not both:\\
\((i)\) The rectangle \(R\) contains an open oriented left right crossing.\\
\((ii)\) The rectangle \(R\) contains a closed dual oriented top bottom crossing.
\end{Theorem}

The paper is organized as follows: In Section~\ref{pf1}, we prove the result Theorem~\ref{thm1} regarding merging of two cycles using a piecewise merging algorithm. In Section~\ref{pf2}, we prove Theorems~\ref{thm3} and~\ref{thm2} regarding the outermost boundary of star and plus connected components, respectively. In Section~\ref{pf6}, we prove Theorem~\ref{thm6} regarding the mutual exclusivity of left right and top bottom crossings in unoriented bond percolation. Finally, in Sections~\ref{pf7} and~\ref{pf8}, we prove the corresponding results for oriented bond percolation. In Section~\ref{pf7}, we obtain preliminary properties regarding the outermost boundary with orientation needed for proving Theorem~\ref{thm8}. In Section~\ref{pf8}, we then prove Theorem~\ref{thm8}.

\section{Proof of Theorem~\ref{thm1}}\label{pf1}

\emph{Proof of Theorem~\ref{thm1}}: If every edge of the cycle~\(C\) either belongs to~\(D\) or is contained in the interior of~\(D,\) then the desired cycle~\(E = D.\) Similarly, if every edge of the cycle~\(D\) either belongs to \(C\) or is contained in the interior of~\(C,\) then \(E = C.\) In what follows, we suppose that the cycle~\(C\) contains at least one edge in the exterior of \(D\) and similarly, the cycle~\(D\) also contains at least one edge in the exterior of~\(C.\)

To merge the cycles \(C\) and \(D,\) we use bridges. Let \(P \subset C\) be any path of edges contained in the cycle~\(C.\) We say that \(P = B(P,D)\) is a \emph{bridge} for cycle~\(D\) if the endvertices of the path \(P\) belong to~\(D\) and every other vertex in~\(P\) lies in the exterior of~\(D.\) In particular, every edge of \(P\) lies in the exterior of the cycle~\(D.\)

We start with cycles \(F_1 = D = (e_1,\ldots,e_s)\) and \(C = (f_1,\ldots,f_t).\) In the first step, we identify a bridge \(P_1\) for the cycle~\(F_1\) contained in the cycle~\(C.\) In the second step, we merge the cycle \(F_1\) with the bridge \(P_1\) to get a new cycle~\(F_2.\) We then repeat the above procedure with the cycle \(F_2\) and continue this process iteratively until all the edges of the cycle \(C\) exterior to the cycle~\(D\) are exhausted. The final cycle obtained is the desired cycle~\(E.\)\\\\
\emph{\underline{Step 1}: Extracting the bridge \(P_1\) from the cycle \(C\)}\\
For \(1 \leq j \leq t-1,\) let \(v_{j}\) be the endvertex common to the edges \(f_j\) and \(f_{j+1}\) and let \(v_0\) be the endvertex common to edges~\(f_1\) and~\(f_t.\) Thus the edge~\(f_j\) in the cycle \(C\) has endvertices \(v_{j-1}\) and \(v_j\) for \(1 \leq j \leq t.\)


Let \(j_1 \) be the least index \(j \geq 1\) so that the edge~\(f_j \in C\) lies in the exterior of the cycle \(D.\) We then have that one endvertex~\(v_{j_1-1}\) of~\(f_{j_1}\) belongs to \(D\) and the other endvertex \(v_{j_1}\) lies in the exterior of~\(D.\) Without loss of  generality, we assume that \(j_1 = 1.\)

\begin{figure}[tbp]
\centering
\includegraphics[width=3.5in, trim= 50 400 100 175, clip=true]{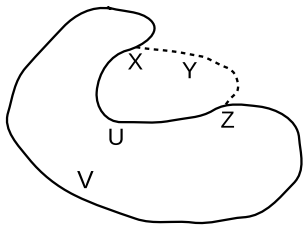}
\caption{Merging the cycle \(XUZVX\) with the segment \(XYZ.\)}
\label{cyc_fig}
\end{figure}

If every edge \(f_k, k \geq 2,\) has both its endvertices in the exterior of the cycle~\(F_1,\) then all the edges of the cycle \(C\) lie in the exterior the cycle \(F_1.\) In particular, since the cycle \(D\) is contained in the interior of the cycle \(F_1,\) we have that every vertex of the cycle \(C\) apart from the vertex \(v_0\) lies in the exterior of the cycle \(D.\) This is a contradiction since we assume that~\(C\) and~\(D\) have more than one vertex in common. Therefore there exists an edge~\(f_k, k \geq 2\) such that the endvertex~\(v_k \neq v_0\) of the edge~\(f_k\) belongs to~\(D.\) Let~\(f_{k_1}\) be the edge with the least such index. We then have that the endvertex \(v_{k_1-1}\) of~\(f_{k_1}\) lies in the exterior of~\(D.\) The path of edges \(P_1 = (f_1,f_2,\ldots,f_{k_1}) \subset C\) is therefore a bridge for the cycle~\(D\) with endvertices \(v_0\) and \(v_{k_1}.\)

In Figure~\ref{cyc_fig}, the cycle \(F_1\) is represented as~\(XUZVX.\) The bridge \(P_1 = XYZ\) with~\(X\) representing the endvertex~\(v_0\) and~\(Z\) representing the endvertex~\(v_{k_1}.\)\\\\
\emph{\underline{Step 2}: Merging the bridge \(P_1\) with the cycle \(F_1\)}\\
The bridge \(P_1\) has endvertices \(v_0\) and \(v_{k_1}.\) Both \(v_0\) and \(v_{k_1}\) also belong to the cycle \(F_1.\) Let \(F_1 = Q_1 \cup R_1\) be the union of two paths where both~\(Q_1\) and~\(R_1\) have endvertices~\(v_0\) and~\(v_{k_1}.\) Let \(G_1=P_1 \cup Q_1\) and \(H_1 = P_1 \cup R_1\) be the two cycles obtained by the union of the bridge \(P_1\) with the two subpaths~\(Q_1\) and~\(R_1\) of the cycle~\(D.\) Exactly one of the cycles, say  \(H_1,\) contains the cycle~\(D\) in the interior and the cycles~\(G_1\) and~\(D\) have mutually disjoint interiors.

In Figure~\ref{cyc_fig}, the paths \(Q_1\) and \(R_1\) are respectively represented by the segments \(XUZ\) and \(XVZ.\) The union of the paths \(P_1 \cup R_1\) is the cycle \(H_1\) which contains the cycle~\(D\) in its interior. We define \(H_1\) to be the cycle obtained at the end of the first iteration. We have the following properties regarding the cycle \(H_1.\)\\
\((a1)\) The cycle \(H_1\) contains only edges from the cycles~\(C\) and \(D.\)\\
\((a2)\) Every edge of the cycle \(D\) either belongs to \(H_1\) or is contained in the interior of~\(H_1.\) Therefore the interior of the cycle \(D\) is contained in the interior of the cycle \(H_1.\)\\
\((a3)\) The cycle \(H_1\) contains at least one edge of \(D\) lying in the exterior of the cycle \(C.\)\\
\emph{Proof of \((a1)-(a3)\) for the cycle \(H_1\)}: The properties \((a1)-(a2)\) are true by construction and the property \((a3)\) is true since every edge in the bridge \(\emptyset \neq P_1 \subset H_1\) lies in the exterior of the cycle \(D.\) \(\qed\)

To proceed to the next step of the iteration, we set \(F_2 = H_1\) and repeat the above procedure with \(F_1\) replaced by~\(F_2.\) Again the cycle \(H_2\) obtained at the end of the iteration step satisfies \((a1)-(a3).\) This procedure continues for a finite number of steps until we obtain a final cycle~\(H_n.\)

It remains to see that the cycle~\(H_n\) is the desired cycle~\(E\) mentioned in the statement of the Theorem. Since~\(H_n\) satisfies property~\((a1),\) we have that the cycle~\(H_n\) contains only edges from the cycles~\(C\) and~\(D.\) Since~\(H_n\) also satisfies property~\((a3),\) we have the property~\((ii)\) in the statement of the theorem is true.

To see that \((i)\) is true, we argue as follows. By definition, the interior of the cycle~\(D\) is contained in the interior of the cycle~\(H_n.\) If there exists an edge of the cycle~\(C\) lying in the exterior of the cycle~\(H_n,\) then we could extract another bridge from cycle \(C\) and the procedure above would not have terminated. Thus every edge in the cycle \(C\) either belongs to \(H_n\) or is contained in the interior of the cycle \(H_n.\) This proves property~\((i).\)

To see the uniqueness of the cycle~\(H_n,\) suppose that there is another cycle \(K \neq H_n\) that satisfies the statement of the Theorem. Without loss of generality, the cycle \(K\) contains an edge \(e\) in the exterior of \(H_n.\) The edge \(e \in C\cup D\) and suppose \(e \in C.\) This means that at least one edge of \(C\) lies in the exterior of \(H_n,\) a contradiction since \(H_n\) satisfies property \((i).\)~\(\qed\)

\section{Proof of Theorems~\ref{thm3} and~\ref{thm2}} \label{pf2}
We first prove Theorem~\ref{thm3} and obtain Theorem~\ref{thm2} as a Corollary. The first step in the proof of Theorem~\ref{thm3} is to obtain large cycles surrounding each occupied square in \(C(0).\) We recall that \(G_0\) is the graph with vertex set being the set of all corners of the squares \(\{S_k\}_{k \geq 0}\) in \(C(0)\) and edge set consisting of the edges of the squares \(\{S_k\}_{k \geq 0}\) in~\(C(0).\) Also \(\{J_k\}_{1 \leq k \leq M} \subset \{S_j\}\) denotes the set of occupied squares belonging to \(C(0).\) We have the following Lemma.

\begin{Lemma}\label{outer} For every \(J_k \in C(0), 1 \leq k \leq M,\) there exists a unique cycle~\(D_k\) in~\(G_0\) satisfying the following properties.\\
\((a)\) The square \(J_k\) is contained in the interior of \(D_k.\)\\
\((b)\) Every edge in the cycle \(D_k\) is a boundary edge adjacent to one occupied square of \(C(0)\) in the interior and one vacant square in the exterior.\\
\((c)\) If \(C\) is any cycle in \(G_0\) that contains \(J_k\) in the interior, then every edge in~\(C\) either belongs to~\(D_k\) or is contained in the interior.\\
\end{Lemma}
We prove in Theorem~\ref{thm6} that every edge of~\(D_k\) is also an outermost boundary edge in the graph~\(G_0.\) We therefore denote~\(D_k\) to be the \emph{outermost boundary cycle} containing the square \(J_k \in C(0).\)

\emph{Proof of Lemma~\ref{outer}}: Fix \(1 \leq k \leq M.\) We first construct a large cycle~\(C_{fin}\) by merging together all cycles containing the square~\(J_k\) in their interior. We then show that the cycle \(C_{fin}\) satisfies properties \((a)-(c).\)

\emph{Construction of the cycle~\(C_{fin}\)}: Let \({\cal E}\) be the set of all cycles in the graph~\(G_0\) satisfying property~\((a);\) i.e., if \(C\) is a cycle containing the square~\(J_k\) in its interior then \(C \in {\cal E}.\) The set~\({\cal E}\) is not empty since the cycle~\(F_0\) formed by the four edges of~\(J_k\) belongs to~\({\cal E}.\) We merge cycles in~\({\cal E}\) one by one using Theorem~\ref{thm1} to obtain a final cycle~\(C_{fin}.\)


We set \({\cal E}_0 = {\cal E}\) and pick a cycle~\(H_0\) in~\({\cal E}_0\setminus F_0\) using a fixed procedure; for example, using an analogous iterative procedure as described in Section~1 of Ganesan (2014) for choosing paths. We merge cycles \(F_0\) and \(H_0\) and from Theorem~\ref{thm1}, we obtain a cycle \(F_1\) in~\(G_0\) containing both~\(F_0\) and~\(H_0\) in its interior. The cycle \(F_1\) also contains the square \(S_0\) in its interior and so satisfies property~\((a).\) If it also satisfies property~\((c),\) then terminate the procedure and output~\(F_1.\)

If the cycle \(F_1\) does not satisfy property \((c),\) then there exists a cycle in~\({\cal E}'_1 = {\cal E}_0 \setminus \{F_0, H_0, F_1\}\) satisfying property~\((a)\) but containing at least one edge in the exterior of the cycle~\(F_1.\) Let \({\cal E}_1 \subset {\cal E}'_1\) be the set of all such cycles and pick one such cycle \(H_1\) using a fixed procedure. Merger \(F_1\) and~\(H_1\) using Theorem~\ref{thm1} to get a new cycle~\(F_2.\)

Repeat the procedure above with the cycle \(F_2\) and this procedure proceeds only for a finite number of steps, since the sets \(\{{\cal E}_j\}\) form a strictly decreasing sequence of subsets of~\({\cal E}\) and the set~\({\cal E}\) has finite number of elements. Let~\(C_{fin}\) be the final cycle obtained at the end of the procedure. It remains to see that the cycle \(C_{fin}\) satisfies properties \((a)-(c).\)

\emph{Proof of properties \((a)-(c)\)}: By construction, the cycle \(C_{fin}\) obtained above is unique and satisfies property~\((a).\) It also satisfies property~\((c)\) because if it did not, the procedure above would not have terminated. It only remains to see that property~\((b)\) is true.

Suppose there exists an edge~\(e\) of the cycle~\(C_{fin}\) that is not a boundary edge. Since~\(e\) belongs to the graph~\(G_0,\) the edge~\(e\) is adjacent to an occupied square \(A_1 \in \{J_i\}.\) But since \(e\) is not a boundary edge, the other square \(A_2 \in \{J_i\}\setminus A_1\) containing~\(e\) as an edge is also occupied. One of these squares, say \(A_1,\) is contained in the interior of \(C_{fin}\) and the other square~\(A_2,\) is contained in the exterior.

The square~\(A_2\) and the cycle~\(C_{fin}\) have the edge \(e\) in common and thus more than one vertex in common. We use Theorem~\ref{thm1} to obtain a larger cycle~\(C_{lar}\) containing both~\(C_{fin}\) and~\(A_2\) in the interior. Since~\(A_2\) contains at least one edge in the exterior of~\(C_{fin},\) the cycle~\(C_{lar}\) also contains at least one edge in the exterior of~\(C_{fin}.\) Moreover, the cycle~\(C_{lar}\) also contains the occupied square~\(J_k\) in its interior and therefore satisfies property~\((a).\) But since~\(C_{fin}\) satisfies property~\((c),\) this is a contradiction. Thus every edge~\(e\) of~\(C_{fin}\) is a boundary edge.

By the same argument above, we also see that the edge~\(e\) cannot be adjacent to an occupied square in the exterior of the cycle~\(C_{fin}.\) Thus~\(e\) is adjacent to an occupied square in the interior of~\(C_{fin}\) and a vacant square in the exterior. Thus the cycle~\(C_{fin}\) satisfies property~\((b).\)\(\qed\) 

\emph{Proof of Theorem~\ref{thm3}}: We claim that the set of distinct cycles in the set~\({\cal D} := \cup_{J_k \in C(0)} \{D_k\}\) obtained in Lemma~\ref{outer} is the desired outermost boundary~\(\partial_0\) and satisfies the properties~\((i)-(v)\) mentioned in the statement of the theorem.

The properties \((iv)\) and \((v)\) follow from Lemma~\ref{outer}. To see that \((iv)\) is satisfied, let \(J_k \in C(0)\) be any occupied square. The outermost boundary cycle \(D_k\) satisfies property \((a)\) of Lemma~\ref{outer} and so contains the square~\(J_k\) in its interior. This proves that~\((iv)\) is true. To prove \((v),\) let \(e \in D_k\) be any edge. Since the cycle \(D_k\) satisfies property~\((b)\) of Lemma~\ref{outer}, the edge \(e\) satisfies~\((v).\)

In what follows, we prove \((iii), (ii)\) and \((i)\) in that order.\\
\emph{Proof of \((iii)\)}: Consider two cycles \(D_{k_1} \neq D_{k_2}.\)  We first see that the cycles \(D_{k_1}\) and \(D_{k_2}\) have mutually disjoint interiors. We consider various possibilities.\\
\((p1)\) Every edge in the cycle~\(D_{k_2}\) is either belongs to or is contained in the interior of the cycle \(D_{k_1}.\)\\
\((p2)\) Every edge in the cycle~\(D_{k_1}\) is either belongs to or is contained in the interior of the cycle \(D_{k_2}.\)\\
\((p3)\) There are edges \(e_1,e_2 \in D_{k_2}\) such that~\(e_1\) lies in the interior of cycle~\(D_{k_1}\) and~\(e_2\) lies in the exterior of~\(D_{k_1}.\)\\
\((p4)\) There are edges \(f_1,f_2 \in D_{k_1}\) such that \(f_1\) lies in the interior of cycle \(D_{k_2}\) and\(f_2\) lies in the exterior of~\(D_{k_2}.\)\\
If none of the above possibilities hold, then the cycles \(D_{k_1}\) and \(D_{k_2}\) have mutually disjoint interiors.

To eliminate possibilities \((p1)-(p2)\) we argue as follows. Suppose \((p1)\) holds. The cycle \(D_{k_1}\) then contains the square~\(J_{k_2} \in C(0)\) is its interior and since \(D_{k_1} \neq D_{k_2},\) the cycle~\(D_{k_1}\) also contains an edge lying in the exterior of~\(D_{k_2}.\) This contradicts the fact that \(D_{k_2}\) satisfies property \((c)\) of Lemma~\ref{outer}. This eliminates possibility~\((p1)\) and an analogous argument holds for \((p2).\)


We eliminate possibilities \((p3)-(p4)\) as follows. If \((p3)\) holds, then the edge \(e_2\) belongs to a path \(P_2 \subset D_{k_2}\) whose every edge lies in the exterior of the cycle \(D_{k_1}.\)  The path \(P_2 \neq D_{k_2}\) since there is at least one edge \(e_1 \in D_{k_2}\) lying in the interior of the cycle \(D_{k_1}.\) From Theorem~\ref{thm1}, we then obtain a cycle~\(E_{12}\) containing both the cycles~\(D_{k_1}\) and~\(D_{k_2}\) in its interior. The cycle \(E_{12}\) also contains an edge \(e_{12} \in D_{k_2}\) lying in the exterior of the cycle~\(D_{k_1}.\) Moreover, the cycle~\(E_{12}\) contains the occupied square \(J_{k_1} \in C(0)\) in its interior and so satisfies property~\((a)\) in Lemma~\ref{outer}. But this is a contradiction since the cycle~\(D_{k_1}\) satisfies property~\((c)\) of Lemma~\ref{outer}. This eliminates possibility~\((p3)\) and an analogous argument holds for~\((p4).\)

We have obtained that \(D_{k_1}\) and \(D_{k_2}\) have mutually disjoint interiors. If they share more than one vertex in common, we again merge them as described in the previous paragraph and obtain a contradiction. Thus the cycles \(D_{k_1} \neq D_{k_2}\) have at most one vertex in common and have mutually disjoint interiors and this proves~\((iii).\)~\(\qed\)

\emph{Proof of \((ii)\)}: We use the fact that the graph~\(G_0\) is connected. To see this is true, let \(u_1\) and \(u_2\) be vertices in \(G_0.\) Each \(u_i, i = 1,2\) is a corner of an occupied square \(S_i \in C(0)\) and by definition, there is a star connected \(S-\)path of squares connecting \(S_1\) and \(S_2,\) consisting only of squares in~\(C(0).\) Thus there exists a path in \(G_0\) from \(u_1\) to \(u_2.\)

To see that \({\cal D} = \cup_{S_k \in C(0)} \{D_k\}\) is a connected subgraph of the graph~\(G_0,\) we let~\(v_1\) and~\(v_2\) be vertices in~\({\cal D}\) that belong to cycles \(D_{r_1}\) and \(D_{r_2},\) respectively, for some \(r_1\) and~\(r_2.\) If \(r_1 = r_2,\) then \(v_1\) and \(v_2\) are connected by a path within the cycle~\(D_{r_1}.\) If \(r_1 \neq r_2,\) let \(P_{12} = (q_1,q_2,\ldots,q_{t-1},q_{t} )\) be a path of edges in~\(G_0\) with endvertices \(v_1\) and~\(v_2.\) We iteratively construct a path~\(Q_{12}\) from~\(P_{12}\) using only edges of cycles in~\({\cal D}.\) Every edge \(q_i\) in~\(P_{12}\) is the edge of an occupied square of \(C(0)\) and so either belongs to a cycle in~\({\cal D}\) or is contained in the interior of some cycle in~\({\cal D}.\) Without loss of generality, we assume that the first edge \(q_1 \in P_{12}\) either belongs to the cycle~\(D_{r_1}\) or is contained in the interior of~\(D_{r_1}.\)

In the first step of the iteration, we let \(s_0 = r_1,i_0=1\) and let~\(i_1 \) be the first time the path~\(P_{12}\) leaves the cycle~\(D_{s_0};\) i.e., let \[i_1 = \min\{i \geq i_0+1 : q_{i} \text{ belongs to exterior of } D_{s_0}\}.\] The edge \(q_{i_1}\) has one endvertex \(w_{i_1}\) in~\(D_{s_0}\) and  the other endvertex lies in the exterior of~\(D_{s_0}.\) Let \(T_1 \subset D_{s_0}\) be a path consisting only of edges in the cycle~\(D_{s_0},\) with endvertices \(v_1\) and \(w_{i_1}.\) This completes the first step of the iteration.

For the second step, we use the fact that the edge \(q_{i_1+1}\) of the path~\(P_{12}\) lies in the exterior of the cycle~\(D_{s_0}\) but belongs to some occupied square of the component~\(C(0).\) We therefore have that either \(q_{i_1+1}\) belongs to some cycle~\(D_{s_1}\) or is contained in its interior. Also the cycles \(D_{s_1}\) and~\(D_{s_0}\) meet at~\(w_{i_1}\) and have no other vertex in common.

Repeating the same procedure above, let \[i_2 = \min\{i \geq i_1+1 : q_{i+1} \text{ belongs to exterior of } D_{s_1}\}\] be the first time \(P_{12}\) leaves the cycle~\(D_{s_1}\) and obtain a path \(T_2 \subset D_{s_1}\) with endvertices~\(w_{i_1}\) and~\(w_{i_2}.\) As before \(w_{i_2}\) is the endvertex of the edge \(q_{i_2}\) belonging to the cycle \(D_{s_1}.\) We continue the above procedure for a finite number of steps \(m,\) until we reach~\(v_2.\) By construction, the path \(T_i\) obtained at step~\(i, 2 \leq i \leq m\) is connected to \(\cup_{1 \leq j \leq i-1} T_j.\) The final union of paths \(\cup_{1 \leq i \leq m}T_{i}\) is therefore a connected graph containing only edges in \({\cal D}\) and also containing the vertices~\(v_1\) and~\(v_2.\) This proves~\((ii).\) \(\qed\)




\emph{Proof of \((i)\)}: We first show that every edge in the union of the cycles \({\cal D} = \cup_{J_k \in C(0)} \{D_k\}\) is an outermost boundary edge. If \(e\) is an edge of a cycle \(D_k \in {\cal D}\) we have that \(e\) is an edge of an occupied square \(J_e \in C(0)\) contained in the interior of~\(D_k\) and is also an edge of a vacant square \(W_e\) in the exterior of \(D_k.\) If \(D_e \in {\cal D}\) denotes the outermost boundary cycle containing the square~\(J_e,\) then from~\((iii)\) above we must have that \(D_e = D_k.\) This is because if \(D_e \neq D_k,\) then \(D_e\) and \(D_k\) have mutually disjoint interiors. This cannot happen since both \(D_e\) and \(D_k\) contain the square \(J_e\) in the interior.

If there exists a cycle \(C\) in the graph~\(G_0\) that contains the edge~\(e\) in the interior, then both the squares \(J_e\) and \(W_e\) are contained in the interior of~\(C.\) Since~\(W_e\) is exterior to the cycle~\(D_e = D_k,\) the cycle~\(C\) contains at least one edge in the exterior of \(D_e.\) This contradicts the fact that the cycle \(D_e\) satisfies property~\((c)\) of Lemma~\ref{outer}. Thus \(e\) is an outermost boundary edge.

We now argue that no other edge apart from edges of cycles in \({\cal D}\) can belong to the outermost boundary. Suppose \(e_1 \notin {\cal D}\) is an edge of the graph~\(G_0\) belonging to some occupied square \(A_1\in C(0).\) From property~\((iv),\) the square~\(A_1\) is contained in the interior of some cycle \(D_k \in {\cal D}.\) Since \(e_1 \notin {\cal D},\) the edge~\(e_1\) does not belong to~\(D_k\) and is necessarily contained in the interior of~\(D_k.\) This proves \((i).\) \(\qed\)


\subsubsection*{Circuit \(C_{out}\) containing the outermost boundary \(\partial_0\)}
To obtain the circuit \(C_{out},\) we first compute the cycle graph \(H_{cyc}\) as follows. Let \(E_1,E_2,...,E_n\) be the distinct outermost boundary cycles in~\({\cal D}.\) Represent~\(E_i\) by a vertex~\(i\) in~\(H_{cyc}.\) If~\(E_i\) and~\(E_j\) share a corner, we draw an edge~\(e(i,j)\) between \(i\) and~\(j.\) We have the following properties regarding the graph \(H_{cyc}.\)\\
\((y1)\) Let \(P = (e(i_1,i_2),e(i_2,i_3),...,e(i_{m-1},i_m))\) be a path of edges in the graph~\(H_{cyc},\) where \(i_j, 1 \leq j \leq m\) are vertices in \(H_{cyc}.\) Let \(u \in E_{i_1}\) and \(v \in E_{i_m}\) be any two vertices. There is a path \(P_{uv}\) with endvertices \(u\) and \(v\) and consisting only of edges of the cycles \(\{E_{i_k}\}_{1 \leq k \leq m}.\)\\
\((y2)\) The graph \(H_{cyc}\) is acyclic and connected and we record the following for future use.
\begin{equation}~\label{hcyc}
\text{We have that the cycle graph } H_{cyc} \text{ is a tree.}
\end{equation}

\emph{Proof of~\((y1)-(y2)\)}:
To obtain the path \(P_{uv},\) we proceed as follows. Set \(w_1 = u, w_{m+1} = v\) and for \(2 \leq j \leq m,\) let \(w_j\) denote the vertex common the cycles~\(E_{i_{j-1}}\) and~\(E_{i_j}.\) For \(1 \leq j \leq m,\) the vertices~\(w_{j}\) and~\(w_{j+1}\) both belong to the cycle~\(E_{i_j}\) and so there is a path~\(Q_{j}\) with endvertices~\(w_{j}\) and~\(w_{j+1}\) consisting only of edges of the cycle~\(E_{i_j}.\) Since the cycles \(\{E_i\}\) are edge disjoint (property \((iii)\)), the paths \(\{Q_j\}_{1 \leq j \leq m}\) are edge disjoint. Therefore the union of the paths~\(\cup_{j=1}^{m}Q_{j}\) is a path with endvertices \(w_1\) and \(w_{m+1}\) and containing only edges in the cycles~\(\{E_{i_j}\}_{1 \leq j \leq m}.\)

To prove \((y2),\) we use property \((ii)\) of the outermost boundary to obtain that the union of the cycles \(\cup_{1 \leq i \leq n} E_i\) is a connected graph. Therefore, the graph \(H_{cyc}\) is connected. We use property \((y1)\) to prove that \(H_{cyc}\) is acyclic.


\begin{figure}[tbp]
\centering
\includegraphics[width=3.5in, trim= 100 350 100 200, , clip=true]{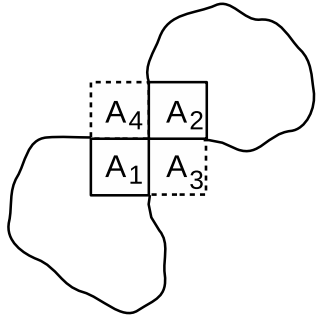}
\caption{Only two cycles of the outermost boundary can meet at a single point.}
\label{two_cyc_meet}
\end{figure}


Suppose \(H_{cyc}\) contains a cycle \(C = (e(r_1, r_2),...,e(r_s, r_1)).\) Again \(1 \leq r_i \leq n, i = 1,2,\ldots,s\) are vertices of the graph~\(H_{cyc}.\) Using property~\((iii)\) of the outermost boundary, we have that the cycles \(E_{r_1} = (u_1,u_2,\ldots,u_m,u_1)\) and~\(E_{r_2}\) have exactly one vertex, say \(u_{j_2},\) in common. Similarly, the cycles~\(E_{r_1}\) and~\(E_{r_s}\) have exactly one vertex,~\(u_{j_s},\) in common.

We have that the indices \(j_2\neq j_s\) since three boundary cycles cannot meet at a point. We assume that some cycle \(E_j, j \neq r_1,r_2\) also contains the vertex~\(u_{j_2}\) and obtain a contradiction as follows. All the four edges~\(g_i, 1\leq i \leq 4,\) containing \(u_{j_2}\) as an endvertex belongs to either \(E_{r_1}\) or \(E_{r_2}.\) The illustrated in Figure~\ref{two_cyc_meet} where the occupied square \(A_1 \in \{S_k\}\) is in the interior of the cycle~\(E_{r_1}\) and the occupied square~\(A_2 \in \{S_k\}\) is in the interior of~\(E_{r_2}.\) The vertex~\(u_{j_2}\) is the vertex common to all the four squares~\(\{A_i\}_{1 \leq i \leq 4}.\) Therefore if the cycle~\(E_j\) contains \(u_{j_2}\) as an endvertex, then \(E_j\) contains at least one of the edges \(g_i, 1 \leq i \leq 4.\) But this contradicts the property~\((iii)\) that the cycles in the outermost boundary have at most one vertex in common.



Let \(Q_1\) and \(R_1\) be the two subpaths of \(E_{r_1}\) with endvertices~\(u_{j_2}\) and~\(u_{j_s}\) so that \(Q_1 \cup R_1 = E_{r_1}.\) The vertex \(u_{j_2} \in E_{r_2}\) and \(u_{j_s} \in E_{r_s}\) and so by property~\((y1)\) above, there exists a path~\(P_{2s}\) with endvertices~\(u_{j_2}\) and~\(u_{j_s},\) consisting only of edges in \(\{E_{r_i}\}_{2 \leq i \leq s}.\) From the property~\((iii)\) of the outermost boundary, we have that the cycles~\(\{E_i\}\) have mutually disjoint interiors. Therefore every edge in the cycle \(E_{r_i}, 2 \leq i \leq s\) lies in the exterior of the cycle~\(E_{r_1}.\) In particular, every edge in the path \(P_{2s}\) lies in the exterior of the cycle~\(E_{r_1}.\)

In Figure~\ref{cyc_fig}, we illustrate the cycle~\(E_{r_1}\) as \(XUZVX\) and the paths~\(R_1\) and~\(Q_1\) are respectively denoted by the segments~\(XVZ\) and~\(XUZ.\) The path~\(P_{2s}\) lying in the exterior of the cycle~\(E_{r_1}\) is the segment~\(XYZ.\) One of the cycles~\(P_{2s} \cup Q_1\) or~\(P_{2s} \cup R_1\) contains the cycle~\(E_{r_1}\) in its interior. We call this cycle~\(C_{12}.\) In Figure~\ref{cyc_fig}, the cycle \(C_{12} = P_{2s} \cup R_1\) is represented by~\(XYZVX\) and contains the cycle~\(E_{r_1} = XUZVX\) in its interior.

Let~\(J_{k} \in \{S_j\}\) be any occupied square of the component~\(C(0)\) in the interior of the cycle~\(E_{r_1}.\) We have from Lemma~\ref{outer} that the cycle~\(E_{r_1} = D_k,\) the outermost boundary cycle containing the square~\(J_k\) and satisfies properties~\((a), (b)\) and \((c)\) mentioned in the statement of Lemma~\ref{outer}. The cycle~\(C_{12}\) also contains~\(J_k\) in its interior and thus satisfies property~\((a).\) Moreover, it contains at least one edge in the exterior of the cycle~\(E_{r_1}\) contradicting the fact that~\(E_{r_1}\) satisfies property~\((c).\) Thus the graph~\(H_{cyc}\) is acyclic. This proves \((y2).\)~\(\qed\)

Using (\ref{hcyc}), we obtain the desired circuit for the outermost boundary~\(\partial_0,\) iteratively, by considering an increasing sequence of tree subgraphs of the tree~\(H_{cyc}.\) The vertex set of \(H_{cyc}\) is \(\{1,2,\ldots,n\}\) and for vertex \(v \in \{1,2,\ldots,n\},\) let~\({\cal N}(v)\) be the neighbours of \(v\) in the tree~\(H_{cyc}.\) Set \(q_1 = 1,\) \(H_1 = \{q_1\}\) and for \(i \geq 1,\) let \(V(H_{i})\) is the vertex set of the graph~\(H_{i}.\) For \(i \geq 1,\) we have the following properties regarding the graph \(H_i.\)\\
\((z1)\) We have that \(H_i\) is a tree subgraph of \(H_{cyc}\) with~\(V(H_i) = \{q_1,\ldots,q_i\}.\) \\
If \(i \leq n-1,\) the following additional condition holds.\\
\((z2)\) There exists a vertex  \(q_{i+1} \notin V(H_i)\) that is adjacent to exactly one vertex \(v_{i+1} \in V(H_i).\) Pick the least such \(q_{i+1}\) and set \(V(H_{i+1}) = V(H_i) \cup \{q_{i+1}\}.\)\\

\emph{Proof of \((z1)-(z2)\) for \(i=1\)}: The proof of \((z1)\) is true by construction. To see \((z2)\) is true, we argue as follows. The vertex set \(V(H_i)\) of the graph~\(H_i\) satisfies \(\#V(H_i) \leq n-1\) and so there is at least one vertex in \(\{1,2,\ldots,n\}\setminus H_i.\) We recall that \({\cal N}(v)\) denotes the neighbours of the vertex \(v\) in the graph~\(H_{cyc}.\) If \({\cal N}(v) \subset H_i\) for all \(v \in H_i,\) then the graph~\(H_i\) is a (connected) component of the graph \(H_{cyc}\) with \(\#V(H_i) \leq n-1.\) This means that \(H_{cyc}\) is not connected, a contradiction since \(H_{cyc}\) is a tree, (see property \((y2)\) above). Thus there exists at least one vertex \(v_{i+1}\in H_i\) containing a neighbour in \(\{1,2,\ldots,n\}\setminus V(H_i).\) Pick the least such~\(v_{i+1}\) and let~\(q_{i+1}\) be the least indexed neighbour of~\(v_{i+1}\) in \(\{1,2,\ldots,n\}\setminus V(H_i).\)

Suppose now that there is another vertex \(w_{i+1} \in H_i, w_{i+1} \neq v_{i+1}\) that is also adjacent to \(q_{i+1} \notin H_i.\) The graph \(H_i\) is a tree and so the vertices~\(w_{i+1}\) and~\(v_{i+1}\) are connected by a path~\(P_{vw}\) consisting only of edges in~\(H_i.\) Let \(f\) be the edge between \(q_{i+1}\) and \(v_{i+1}\) in the tree~\(H_{cyc}\) and let \(g\) be the edge between the vertices~\(q_{i+1}\) and~\(w_{i+1}\) in the tree~\(H_{cyc}.\) The union \(P \cup \{f,g\}\) contains a cycle consisting of edges in \(H_{cyc},\) a contradiction to~(\ref{hcyc}). \(\qed\)

The property \((z2)\) is used to proceed to the next step of the iteration and the graph \(H_2\) obtained at the end of the second iteration again satisfies properties \((z1)-(z2).\) The above procedure continues for \(n\) steps and the final graph \(H_n = H_{cyc}.\)

We use the graphs \(H_i, 1  \leq i \leq n\) to construct the desired circuit for the outermost boundary iteratively. We recall that \(\{E_i\}_{1 \leq i \leq n}\) are the cycles in the outermost boundary \(\partial_0.\) The graph \(H_1\) contains a single vertex \(\{1\}\) and set \(\Pi_1 = E_1\) to be the circuit obtained at the end of the first iteration. For \(1 \leq i \leq n,\) let \(\Pi_{i}\) be the circuit obtained at the end of the \(i^{th}\) iteration. We have the following properties.\\
\((x1)\) The circuit \(\Pi_{i}\) contains all the edges belonging to the cycles \(E_v, v \in V(H_{i}) = \{q_1,\ldots,q_{i}\}.\)\\
\((x2)\) If \(i \leq n-1,\) then the circuit \(\Pi_{i}\) does not contain any edge from \(\bigcup_{j=i+1}^{n} E_{q_j}\) and shares exactly one vertex with the cycle~\(E_{q_{i+1}}.\)\\
\emph{Proof of \((x1)-(x2)\) for \(i = 1\)}: The proof of \((x1)\) is true by construction. To see \((x2)\) is true, we use the fact that  tree \(H_{i}\) has vertex set \(\{q_1,q_2,\ldots,q_{i}\}.\) From property~\((z2)\) above, we have that the vertex \(q_{i+1}\) is adjacent to exactly one vertex \(v_{i+1} \in V(H_i).\)
The corresponding cycle~\(E_{q_{i+1}}\) therefore shares exactly one vertex with the cycle~\(E_{v_{i+1}} \subset \Pi_i\) and does not share a vertex with any other cycle~\(E_v, v \in H_{i}.\) \(\qed\)

We use property \((x2)\) to proceed to the next step of the iteration.  Let  \(\Pi_{i} = (c_1,c_2,...,c_r)\) be the edges of the circuit~\(\Pi_{i}\) traversed in that order and let \(\Pi_{i}\) meet the cycle \(E_{q_{i+1}} = (d_1,d_2,...,d_l)\) at some vertex \(a\) common to the edges~\(c_1\) and~\(d_l.\) We then form the new circuit \(\Pi_{i+1} = (d_1,d_2,...,d_l,c_1,c_2,...,c_r,c_1).\) The circuit \(\Pi_{i+1}\) also satisfies properties \((x1)-(x2)\) and continuing this way iteratively, the final circuit~\(\Pi_n\) contains the edges of all the cycles \(E_v, 1 \leq v \leq n.\) The circuit \(\Pi_n\) is therefore the desired circuit of the outermost boundary~\(\partial_0.\)~\(\qed\)

\emph{Proof of Theorem~\ref{thm2}}: Let \(D_0\) be the outermost boundary cycle containing the square \(S_0\) as in Lemma~\ref{outer}. It satisfies the conditions \((i)\) and \((ii)\) in the statement of the theorem and is unique and thus \(\partial_0^+ = D_0.\) \(qed\)



\section{Proof of Theorem~\ref{thm6}}~\label{pf6}
We first see that an open left right crossing and a closed dual top bottom crossing cannot occur simultaneously. Let \(P_1\) be an open left right crossing of~\(R.\) If there exists a closed top bottom dual crossing~\(P_2\) of \(R,\) then the paths~\(P_1\) and~\(P_2\) intersect and in particular, there is some edge \(e \in P_1\) which is open but its dual edge \(f(e) \in P_2\) is closed. This is a contradiction.


We henceforth assume that \(R\) does not contain an open left right crossing. Let \({\cal C}\) be the set of all vertices in \(R\) that are connected to some vertex in~\(R_{left}\) by an open path of edges. We recall that~\(R_{left}\) denotes the left edge of the rectangle~\(R\) and every vertex \( (0,i) \in R_{left}\) belongs to~\({\cal C}.\) Also since \(R\) has no open left right crossing, no vertex in \({\cal C}\) belongs to the right edge~\(R_{right}\) of~\(R.\)

We also recall that every vertex in the rectangle~\(R\) is the centre of some square in \(\{S_k\}\) contained in the dual graph~\(G.\) Let \(J_i \subset \{S_k\}\) be the square with centre as~\((0,i) \in R_{left}.\) For \(z \in R \setminus R_{left},\) let \(J_z \in \{S_k\}\) be the square with centre~\(z.\) We say that \(J_z\) is occupied if the corresponding vertex \(z \in {\cal C}.\) Else we set \(J_z\) to be vacant. Let \(C_{left}\) be the plus connected component containing the square~\(J_0 = S_0\) with centre as the origin. Every square \(J_i, 1 \leq i \leq n\) belongs to \(C_{left}.\) From Theorem~\ref{thm2}, we have that the outermost boundary~\(\partial^+_{left}\) of \(C_{left}\) is a single cycle.




We extract a closed dual top bottom crossing as a subpath of \(\partial^+_{left}.\) For \(0 \leq i \leq n,\) let \(h_i(t), h_i(r), h_i(b)\) and \(h_i(l)\) denote the top, right, bottom and left edges of the square~\(J_i.\) Also let \(R_{top}\) and \(R_{bottom}\) be the top and bottom edges of the rectangle~\(R.\) We enumerate the following properties of \(\partial^+_{left}\) needed for future use. \\
\((x1)\) The path \(\Pi_1 = (h_n(t),h_n(l),h_{n-1}(l),\ldots,h_0(l),h_0(b))\) is a subpath of~\(\partial^+_{left}.\)\\
Let \(\Pi_2 := \partial^+_{left} \setminus \Pi_1 = (f_1,\ldots,f_r)\) where the edge \(f_1\) shares an endvertex with the top edge \(h_n(t)\) of \(\Pi_1\) and the edge \(f_r\) shares an endvertex with bottom edge \(h_0(b)\) of \(\Pi_1.\)\\
\((x2)\) If \((x,y)\) is an endvertex of some edge in \(\Pi_2,\) then \(\frac{1}{2} \leq x \leq m-\frac{1}{2}\) and \(-\frac{1}{2} \leq y \leq n+\frac{1}{2}.\) In particular, no edge in subpath \(\Pi_2\) intersects the left edge~\(R_{left}\) of \(R\) and no edge in~\(\Pi_2\) intersects the right edge~\(R_{right}\) of \(R.\)\\
\((x3)\) Let \(\Pi = (h_1,\ldots,h_s))\) be any subpath of \(\Pi_2\) with one endvertex of \(h_1\) lying above \(R_{top}\) and one endvertex of \(h_s\) lying below \(R_{bottom}.\) There are indices \(1 \leq k_1 < k_2 \leq s\) such that \(h_{k_1}\) intersects the top edge of \(R\) and the edge \(h_{k_2}\) intersects the bottom edge of \(R.\)\\

Essentially property \((x3)\) implies that the path \(\Pi_2\) intersects the top edge of \(R\) before intersecting the bottom edge.

\emph{Proof of \((x1)-(x3)\)}: We prove \((x1)\) as follows. The squares \(J_i, 1 \leq i \leq n\) are the left most squares in the plus connected component \(C_{left}\) in the sense that the square \(J_i\) has centre \((0,i)\) in the left edge \(R_{left}\) of the rectangle \(R.\) Suppose the left edge \(h_i(l)\) of \(J_i\) does not belong to the outermost boundary cycle~\(\partial^+_{left}.\) Some edge \(e\) in the cycle~\(\partial^+_{left}\) then intersects the line \(y = i\) at~\((x,i),\) where \(x \leq -\frac{3}{2}.\) This is true because the edge \(h_i(l)\) is contained in the line \(x = -\frac{1}{2}\) and every edge of \(\partial^+_{left}\) intersects the line \(y = i\) at \((\frac{k}{2},i)\) for integer \(k \neq 0.\) The edge \(e\) belongs to a square \(J_v\) whose vertex \(v \in {\cal C}.\) The vertex~\(v\) has \(x-\)coordinate at most~\(-1,\) a contradiction since every vertex in~\({\cal C}\) belongs to~\(R.\)

An analogous argument as above using top most and bottom most squares obtains that~\(h_n(t)\) and~\(h_0(b)\) belong to~\(\partial^+_{left}.\) Since every vertex apart from the endvertices have degree two in \(\Pi_1,\) the path \(\Pi_1\) is a subpath of the cycle~\(\partial^+_{left}.\)

We prove \((x2)\) is true as follows. The bounds for \(y-\)coordinates are true since we only consider edges of squares with centre in the rectangle~\(R.\)
Since the subpath~\(\Pi_1 \subset \partial^+_{left},\) we have that no edge in \(\Pi_2\) intersects the left edge of \(R.\) Because if there exists an edge~\(f\) intersecting the left edge of~\(R,\) then \(f\) shares an endvertex~\(v\) with the left edge \(h_i(l)\) of some square~\(J_i, 0 \leq i \leq n-1.\) Since the edges~\(h_i(l)\) and \(h_{i+1}(l)\) both belong to the cycle~\(\partial^+_{left},\) the vertex  \(v\) then has degree three in~\(\partial^+_{left},\) a contradiction. Thus if \((x_1,y_1)\) is the endvertex of some edge in \(\Pi_2,\) we have \(x_1 \geq \frac{1}{2}.\)

Suppose now that some edge \(g\) in \(\Pi_2\) has endvertex intersecting the line \(x = m+\frac{1}{2}.\) The edge \(g\) is the edge of a square \(J_v\) whose centre \(v\) lies in the right edge \(R_{right}\) of \(R\) and also belongs to the cluster \({\cal C}.\) This means that \({\cal C}\) contains a left right crossing of \(R,\) a contradiction. Thus  the \(x-\)coordinate of the endvertex \((x_1,y_1)\) of any edge in \(\Pi_2\) satisfies \(x_1 \leq m-\frac{1}{2}.\)


To prove \((x3),\) we use the fact that the path \(\Pi\) crosses both the lines \(y = n\) and \(y = 0.\) Moreover, it crosses \(y = n\) before crossing \(y = 0;\) i.e., there are edges \(h_{k_1}\) and \(h_{k_2}\) such that \(k_1 < k_2\) and \(h_{k_1}\) crosses the top line \(y = n\) and the edge \(h_{k_2}\) crosses the bottom line \(y = 0.\) From property \((x2),\) we have that the endvertices of \(h_{k_1}\) lie in the cylinder \(\frac{1}{2} \leq x \leq m-\frac{1}{2}\) and so the edge \(h_{k_1}\) intersects the top edge of \(R\) and similarly, the edge \(h_{k_2}\) intersects the bottom edge of~\(R.\)~\(\qed\)

We use the above properties to extract the necessary dual top bottom crossing.
Let \(k_2\) be the least index such that the edge \(f_{k_2} \in \Pi_2\) crosses the line \(y = 0.\) From property \((x2),\) we have that one endvertex~\(v_{k_2}\) of \(f_{k_2}\) lies in the interior of~\(R.\) The other endvertex \(u_{k_2}\) lies in the exterior of \(R.\) Also, either \(f_{k_2-1}\) or \(f_{k_2+1}\) contain the endvertex \(v_{k_2}.\) If~\(f_{k_2+1}\) contains the endvertex~\(v_{k_2},\) then the edge \(f_{k_2-1}\) contains \(u_{k_2}\) as an endvertex and \(u_{k_2}\) lies below the bottom edge \(R_{bottom}\) of \(R.\) From property \((x3),\) we then have that the subpath~\((f_1,\ldots,f_{k_2-1})\) contains an edge that crosses the line \(y = 0,\) a contradiction to the definition of \(k_2.\)

From the discussion in the previous paragraph, we have that the edge~\(f_{k_2-1}\) contains~\(v_{k_2}\) as an endvertex and since the height \(n\) of the rectangle~\(R\) is at least two, the edge~\(f_{k_2-1}\) does not intersect the top edge of~\(R.\) Therefore from property \((x2),\) both the endvertices of the edge~\(f_{k_2}\) lie in the interior of \(R.\) Since an endvertex the first edge \(f_1\) lies above \(y = n,\) the subpath \((f_1,\ldots,f_{k_2-1})\) crosses the line \(y = n\) and so there exists an index \(k_1 < k_2-1\) such that the edge~\(f_{k_1}\) crosses the top line \(y = n.\) Let~\(r_1\) be the largest index less than~\(k_2-1\) such that~\(f_{r_1}\) crosses the top edge~\(R_{top}\) of~\(R.\) Let~\(v_{r_1}\) be the endvertex of~\(f_{r_1}\) belonging in the interior of \(R.\) Arguing as in the previous paragraph, we have that the edge~\(f_{r_1+1}\) contains~\(v_{r_1}\) as an endvertex and both the endvertices of~\(f_{r_1+1}\) lie in the interior of~\(R.\)


Every edge in the subpath \(\Pi_t = (f_{r_1+1},\ldots,f_{k_2-1})\) has both its endvertices in the interior of the rectangle~\(R.\) Therefore the union \((f_{r_1},\Pi_t,f_{k_2})\) is a dual top bottom crossing of \(R.\) \(\qed\)


\section{Outermost boundary in oriented percolation}~\label{pf7}



We recall that \(R_{left}\) denotes the left edge of the rectangle \(R.\) For simplicity, we denote the square \(S''_z\) containing the vertex \(z \in \mathbb{Z}^2\) as the centre, simply as \(S_z.\)  The edges in \(S_z\) are called dual edges and in this subsection, we do not consider orientation in the dual edges. We introduce the corresponding orientation in the next subsection.

Let \(f_{NE}(z)\) denote the unoriented (dual) edge of the square~\(S_z\) with endvertices \((i-1,j)\) and \((i,j+1)\) and let \(f_{SE}(z)\) denote the unoriented edge with endvertices~\((i,j+1)\) and~\((i+1,j)\) with notations \(NE\) and \(SE\) standing for north east and south east, respectively. Similarly we define north west unoriented edge \(f_{NW}(z)\) as the edge with endvertices \((i+1,j)\) and \((i,j-1)\) and the south west unoriented edge \(f_{SW}(z)\) with endvertices \((i,j-1)\) and \((i-1,j).\) 

For illustration we refer to Figure~\ref{or_mod}\((b)\) above, where \(F_1F_2F_3F_4\) represents the square \(S_{A}\) containing the point \(A = (0,1) \in R_{left}\) as the centre. The segments \(F_1F_2,F_2F_3, F_3F_4\) and \(F_4F_1\) respectively denote the dual edges \(f_{SE}(0,1),f_{SW}(0,1),f_{NW}(0,1)\) and \(f_{NE}(0,1).\)


Let \({\cal C}\) denote the collection of all vertices in the open oriented cluster defined as follows. If \(x  = (0,j) \in R_{left}\) then  \(x \in {\cal C}\) if and only if \(j\) is odd. Since the height \(n\) of the rectangle is at least one, we have that \((0,1) \in {\cal C}\) and so \({\cal C} \neq \emptyset.\) If \(x \in R \setminus R_{left},\) then \(x \in {\cal C}\) if and only if there is an open oriented path~\(P_x\) with one endvertex in \({\cal C} \cap R_{left}\) and the other endvertex as \(x.\) We recall that every open oriented path contains one arrow endvertex and one non arrow endvertex and by construction, the non arrow endvertex of \(P_x\) is some \((0,j) \in {\cal C} \cap R_{left}\) and the arrow endvertex is~\(x.\)

For \(z \in \mathbb{Z}^2,\) we say the square~\(S_z\) is \emph{occupied} if either \(z \in {\cal C}\) and \emph{vacant} otherwise. The resulting union \({\cal Q} = \cup_{z \in {\cal C}} S_z\) of occupied squares is a star connected component. From Theorem~\ref{thm3}, the outermost boundary of \({\cal Q}\) is a unique connected union of cycles \(\cup_{i=1}^{h}C_i\) consisting of dual edges in \(\cup_{z \in {\cal C}} S_z\) and satisfying the following properties:\\
\((a1)\) Every vertex \(z \in {\cal C}\) is in the interior of some cycle \(C_i.\)\\
\((a2)\) For any \(1 \leq i \leq h,\) every edge in the cycle \(C_i\) is a boundary edge adjacent to one occupied square contained in the interior of \(C_i\) and one vacant square in the exterior. \\
\((a3)\) The cycles \(\{C_i\}\) have mutually disjoint interiors and for \(i \neq j,\) the cycles~\(C_i\) and~\(C_j\) intersect at most at one point.

The following is the main result we prove in this Section.
\begin{Theorem}\label{oriented_cycle}
With orientation as introduced in Section~\ref{intro}, every cycle \(C_i, 1 \leq i \leq h,\) in the outermost boundary \(\partial_0\) is an oriented dual cycle; i.e. an oriented dual path with coincident endvertices.
\end{Theorem}
We also derive auxiliary properties along the way used in obtaining the dual crossing in the next section.

In the following two subsections, we consider the outermost boundary without orientation and in the final subsection, we introduce orientation and prove Theorem~\ref{oriented_cycle}.

\subsubsection*{Contiguous block property of the outermost boundary}
In the main result of this subsection, we state and prove the contiguous block property for the outermost boundary. We recall that \(R_{left}\) is the left edge of the rectangle \(R.\)\\\\
\((b1)\) If \(z \in {\cal C} \cap R_{left},\) then the dual edges \(f_{NW}(z)\) and \(f_{NE}(z)\) belonging to the square \(S_z\) are consecutive edges in some cycle \(C_i\) of the outermost boundary. Also, both \(f_{NW}(z)\) and \(f_{NE}(z)\) lie in the exterior of every cycle \(C_j, 1 \leq j \leq h, j\neq i.\)\\
\((b2)\) If vertices \((0,j_1),(0,j_2) \in {\cal C} \cap R_{left}, j_1 < j_2\) both belong to the interior of some cycle~\(C_i\) of the outermost boundary, then every \((0,j) \in {\cal C} \cap R_{left}\) with \(j_1 \leq j \leq j_2\) belongs to the interior of~\(C_i.\)\\

The property \((b2)\) is the contiguous block property which says that the set of vertices of \(R_{left}\) lying in the interior of a boundary cycle \(C_i\) forms a contiguous block.\\\\
\emph{Proof of \((b1)-(b2)\)}: We prove \((b1)-(b2)\)  for the cycle \(C_1\) containing the vertex \((0,1) \in {\cal C} \cap R_{left}\) in its interior. Let \(z = (0,j) \in {\cal C} \cap R_{left}\) be any vertex in the interior of the cycle~\(C_1\) of the outermost boundary. We prove the property for~\(f_{NW}(z)\) and an analogous proof holds for \(f_{SW}(z).\) If the edge~\(f_{NW}(z)\) does not belong to the cycle~\(C_1,\) then it lies in the interior of~\(C_1\) and so both the squares containing~\(f_{NW}(z)\) as an edge lie in the interior of~\(C_1.\) Since the squares~\(S_z\) and~\(S_{z_1}\) with centre \(z_1 = (-1,j-1)\) both contain~\(f_{NW}(z)\) as an edge, the vertex~\(z_1\) lies in the interior of the cycle~\(C_1.\) This is a contradiction since every vertex in the interior of \(C_1\) either lies in the interior of the rectangle~\(R\) or belongs to the boundary and the vertex \(z_1\) lies in the exterior of~\(R.\) This proves \((b1).\)

We assume that \((b2)\) is not true and arrive at a contradiction. Suppose there are integers \(j_1< j_2\) such that the following two statements \((a)-(b)\) hold.~\((a)\)~The vertices \((0,j_1), (0,j_2) \in {\cal C} \cap R_{left}\) both belong to the interior of the cycle \(C_1=  (f_1,f_2,\ldots,f_t)\) of the outermost boundary. \((b)\) Every intermediate vertex \((0,j) \in {\cal C} \cap R_{left},  j_1+1 \leq j \leq j_2-1, \) belongs to the exterior of~\(C_1.\)

We recall that~\(R_{left}\) is the left edge of the rectangle~\(R\) and \({\cal C}\) is the set of vertices in the rectangle~\(R\) which are connected by an oriented open path to some vertex in \({\cal C} \cap R_{left} = \{(0,j), 1\leq j \leq n, j \text{ odd}\}.\) We therefore assume that \(j_1+2 < j_2\) and arrive at a contradiction. Since the vertex \((0,j_1) \in {\cal C} \cap R_{left}\) lies in the interior of the cycle~\(C_1,\) we have from property~\((b1)\) that the \(NE\) edge~\(f_{NE}(0,j_1)\)  associated with the vertex \((0,j_1)\) belongs to \(C_1\) and \(f_{i_1} = f_{NE}(0,j_1)\) for some index \(1 \leq i_1 \leq t.\) The edge~\(f_{i_1}\) contains \((0,j_1+1)\) as an endvertex.

\begin{figure}[tbp]
\centering
\includegraphics[width=3.5in, trim= 0 250 0 175, clip=true]{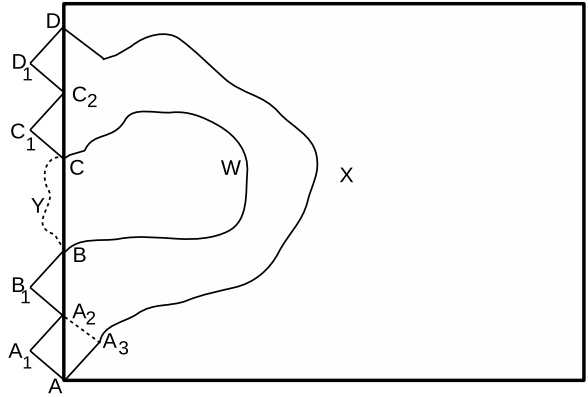}
\caption{Illustration of the square \(S_{0,1} = AA_1A_2A_3\) containing the point \((0,1) \in {\cal C} \cap R_{left}.\) Also depicted is the cycle~\(C_1\) containing the square \(S_{0,1}\) in its interior.}
\label{or_perc11}
\end{figure}

In Figure~\ref{or_perc11}, we illustrate the scenario in the previous paragraph. The points \(A =(0,0),A_2 = (0,2)\) and the mid point of \(AA_2\) is \((0,1).\) The square~\(S_{0,1}\) with centre~\((0,1)\) is denoted by \(AA_1A_2A_3.\) The cycle~\(C_1\) is denoted by the sequence \(AA_1A_2B_1BWCC_1C_2D_1DXA_3A.\) The point \(B = (0,4)\) denotes the vertex~\((0,j_1+1)\) and the edge~\(f_{i_1}\) is denoted by the segment \(BB_1.\)

In an analogous manner, the north west edge \(f_{NW}(0,j_2)\) of the vertex \((0,j_2)\in {\cal C} \cap R_{left}\) also belongs to the cycle~\(C_1\) and moreover \(f_{NW}(0,j_2) = f_{i_2}\) for some index \(1 \leq i_2 \leq r.\) The edge \(f_{i_2}\) contains \((0,j_2-1)\) as an endvertex. In Figure~\ref{or_perc11}, the point \(C\) denotes the vertex~\((0,j_2-1)\) and the edge \(f_{i_2}\) is denoted by the segment \(CC_1.\)

Let \(P_{12}\) be the path containing the union of the edges\\\(\{f_{NW}(0,j),f_{NE}(0,j)\}, j_1+1 \leq j \leq j_2-1.\) From property~\((b1),\) we have that  every edge in~\(P_{12}\) lies in the exterior of the cycle~\(C_1\) and the path~\(P_{12}\) contains \((0,j_1+1)\) and \((0,j_2-1)\) as endvertices. The dotted line \(BYC\) denotes \(P_{12}\) in Figure~\ref{or_perc11}.

The cycle \(C_1= Q_1 \cup R_1\) is the union of two paths with endvertices~\((0,j_1+1)\) and~\((0,j_2-1).\) The union of the paths~\(P_{12} \cup Q_1\) and~\(P_{12} \cup R_1\) are therefore two cycles with the following property. Exactly one of the cycles, say \(P_{12} \cup R_1\) contains the cycle~\(C_1\) in its interior and the other cycle~\(P_{12} \cup Q_1\) has mutually disjoint interior with the cycle~\(C_1.\) In Figure~\ref{or_perc11}, the segment \(BWC\) denotes~\(Q_1\) and the cycles \(P_{12} \cup BWC\)  and \(C_1\) have mutually disjoint interiors.

The cycle \(P_{12} \cup R_1\) contains the cycle~\(C_1\) in its interior and at least one edge in the exterior of \(C_1.\) This contradicts the construction of the cycle~\(C_1,\) which is the outermost boundary cycle containing the occupied square~\(S_{0,1};\) see condition~\((c)\) of Lemma~\ref{outer}.~\(\qed\)


\subsubsection*{Definition and properties of the integers \(\{m_i\}\)}
To study the dual edges with orientation, we need a couple of additional properties regarding the outermost boundary cycles. We use contiguous block property~\((b2)\) and define a increasing sequence of integers
\begin{equation}\label{m_def}
1 = m_1 < m_2 < \ldots < m_h < m_{h+1}
\end{equation}
as follows. We assume that the cycle~\(C_1\) contains the vertex \((0,1) \in {\cal C} \cap R_{left}\) in its interior. We recall that \(R_{left}\) is the left edge of the rectangle~\(R\) and \({\cal C} \cap R_{left}= \{(0,j) : 1 \leq j \leq n, j \text{ odd}\}.\) Let \(m_2-2\) be the largest index \(j\) such that \((0,j)\) lies in the interior of cycle \(C_1.\) The cycle~\(C_1\) contains vertices~\((0,m_1)\) and \((0,m_2-2)\) in its interior and from the contiguous block property \((b2),\) we have that \(C_1\) contains all the vertices \((0,j) \in {\cal C} \cap R_{left}, m_{1} \leq j \leq m_{2}-2\) in its interior and no other vertices of~\(R_{left}\) in its interior.

Since the vertex \((0,m_2-2)\) lies in the interior of the cycle \(C_1,\) we have that \(m_2\) is odd and so the smallest integer \(j\) such that \((0,j)\) lies in the interior of some cycle \(C_j ,2 \leq j \leq h\) is \(j = m_2.\) We assume that the vertex \((0,m_2)\) lies in the  cycle \(C_2\) and proceeding as in the previous paragraph, we obtain an index \(m_3\) such that all vertices \((0,j) \in {\cal C}\cap R_{left}, m_2 \leq j \leq m_3-2\) lie in the interior of the cycle \(C_2.\) Also every other vertex in \({\cal C} \cap R_{left}\) lies in the exterior of~\(C_2.\)

Continuing this way iteratively, we obtain the sequence~\(\{m_i\}.\) The final integer~\(m_{h+1}\) satisfies the following bounds
\begin{equation}\label{m_h_1_def}
n+1 \leq m_{h+1}\leq n+2
\end{equation}
where we recall that~\(n\) is the height of the rectangle \(R.\)\\
\emph{Proof of~(\ref{m_h_1_def})}: To obtain the second inequality in (\ref{m_h_1_def}), we argue as follows. The final cycle~\(C_h\) contains all the vertices \(\{(0,j) \in {\cal C} \cap R_{left}, m_h \leq j \leq m_{h+1}-2\}\) in its interior. The top most vertex in~\(R_{left}\) is~\((0,n)\) and so \(m_{h+1}-2 \leq n.\)

To obtain the lower bound in (\ref{m_h_1_def}), we assume \(m_{h+1} \leq n\) and arrive at a contradiction. Since \(m_{h+1}-2 \leq n-2,\) the vertices \((0,n-1)\) and \((0,n)\) both lie in the exterior of all the cycles \(C_i, 1 \leq i \leq h.\) But one of the integers~\(n-1\) or~\(n\) is odd and so one of the vertices \((0,n-1)\) or \((0,n)\) belongs to \({\cal C} \cap R_{left}.\) Also, by construction, every vertex in \({\cal C}\) lies in the interior of some cycle \(\{C_i\}\) (see property \((a1)\)). This leads to a contradiction.~\(\qed\)

\begin{figure}[tbp]
\centering
\includegraphics[width=3.5in, trim= 50 300 50 100, clip=true]{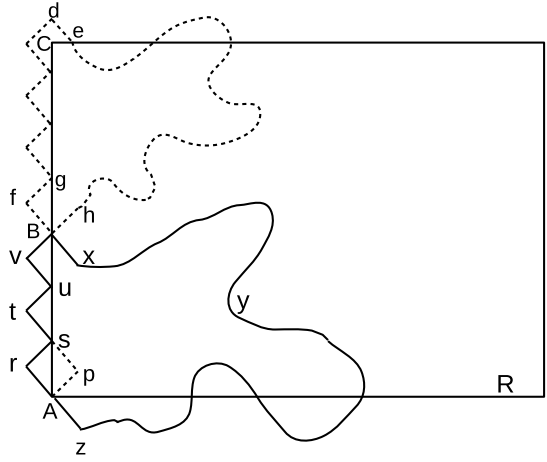}
\caption{The cycles \(C_1\) (thick line) and \(C_2\) (dotted line) are shown meeting at the point \(B\) on the left edge of the rectangle \(R.\) }
\label{c1_c2_unor}
\end{figure}

We have the following properties regarding the integers \(\{m_i\}.\) For each property, we also provide the corresponding illustration in Figure~\ref{c1_c2_unor}.\\
\((d1)\) Fix \(1 \leq i \leq h.\) The path
\begin{eqnarray}
&&\Gamma_i = (f_{NW}(0,m_i),f_{NE}(0,m_i),f_{NW}(0,m_i+1),f_{NE}(0,m_i+1),\nonumber\\
&&\;\;\;\;\;\;\;\;\;\;\;\;\;\;\ldots,f_{NW}(0,m_{i+1}-2),f_{NE}(0,m_{i+1}-2)).\nonumber
\end{eqnarray}
is a subpath of the cycle \(C_i\) with endvertices \((0,m_i-1)\) and \((0,m_{i+1}-1).\) The path \(\Delta_i = C_i \setminus \Gamma_i\) also has endvertices \((0,m_i-1)\) and \((0,m_{i+1}-1).\) Let \(e_i \in \Delta_i\) be the edge sharing an endvertex with the first edge \(f_{NW}(0,m_i) \in \Gamma_i\) and let \(g_i \in \Delta_i\) be the edge sharing an endvertex with the last edge \(f_{NE}(0,m_{i+1}-2) \in \Gamma_i.\) Every edge in \(\Delta_i \setminus \{e_i,g_i\}\) has both endvertices in \(R\setminus R_{left}.\) Every edge in \(\Delta_i\) has both its endvertices to the right of \(R_{left}.\)

The property \((d1)\) is illustrated in Figure~\ref{c1_c2_unor} for the case when there are two cycles~\(C_1\) and~\(C_2\) in the outermost boundary. The cycle \(C_1\) is drawn with thick lines and the cycle \(C_2\) is drawn with dotted lines. The vertices \((0,0),\) \((0,1), (0,2) \) etc are respectively denoted by the points~\(A,\) the mid point of the segment \(As,\)  the point~\(s\)~etc. The segment \(Ar\) represents the first edge \(f_{NW}(0,1) \in C_1\) containing the origin. The path \(ArstuvB\) is the subpath~\(\Gamma_1 \subset C_1.\) The path~\(\Delta_1\) is the wavy segment~\(BxyzA.\) The line segments~\(Bx\) and~\(Az\) respectively denote the edges \(g_1\) and \(e_1.\)  \\\\
\((d2)\) Suppose that \(h \geq 2.\) For \(1 \leq i \leq h-1,\) the cycles \(C_{i}\) and \(C_{i+1}\) intersect only at the vertex \((0,m_{i+1}-1) \in R_{left}\) and have no other vertex in common. For \(2 \leq i \leq h,\) the edge \(e_i = f_{SW}(0,m_i)\) and for \(1 \leq i \leq h-1,\) the edge \(g_i = f_{SE}(0,m_{i+1}-2).\)

In Figure~\ref{c1_c2_unor}, the cycles \(C_1\) and \(C_2\) meet at the point~\(B = (0,m_2-1) = (0,6)\) on the left edge \(R_{left}\) of the rectangle~\(R.\) The line segments \(Bx\) and \(Az\) respectively denote the edges \(g_1 = f_{SE}(0,m_2-2) \in \Delta_1\) and \(e_1 \in \Delta_1.\)  From the figure, we also see that the segment \(Ap\) representing the edge \(f_{SW}(0,1)\) need not necessarily belong to the first cycle~\(C_1.\) Thus the second statement of property~\((d2)\) regarding the edge \(e_i\) need not hold for \(i = 1.\) An analogous argument holds for the edge \(g_h\) belonging to the top most cycle~\(C_h.\)

The edges \(f_{NW}(0,m_2-2)\) and \(e_2 = f_{SW}(0,m_2-2)\) belonging to the cycle~\(C_2\) are represented by the segments~\(Bf\) and~\(Bh.\) The corner~\(C\) of the rectangle~\(R\) belongs to \({\cal C} \cap R_{left}\) and so represents the vertex \((0,m_3-2) = (0,13).\) The edge~\(g_2 = f_{SE}(0,m_3-2)\) represented by the segment~\(de\) lies in the exterior of the rectangle~\(R.\)\\\\
\((d3)\) Suppose \(h \geq 3.\) The cycle \(C_1\) intersects only the cycle \(C_2\) and no other cycle \(C_j,2 \leq j~\leq~h.\) The cycle \(C_h\) intersects only the cycle~\(C_{h-1}\) and no other cycle \(C_j, 1 \leq j \leq h-2.\) For every \(2 \leq i \leq h-1,\) the cycle~\(C_i\) intersects only the cycles~\(C_{i-1}\) and~\(C_{i+1}\) and no other cycle in the outermost boundary. \\






The properties \((d1)-(d3)\) are also used in the next subsection to obtain properties regarding the boundary cycles \(\{C_i\}\) with orientation.\\\\
\emph{Proof of \((d1)-(d3)\)}: We prove the statement of \((d1)\) regarding the subpath \(\Gamma_i\) for \(i = 1.\) An analogous proof holds for general~\(i.\) By definition of the integers~\(\{m_i\},\) the vertices \((0,j) \in {\cal C} \cap R_{left}, 1 \leq j \leq m_2-2\) lie in the interior of the cycle~\(C_1.\) We recall that \({\cal C} \cap R_{left} = \{(0,j)  \in R_{left} : j \text{ odd}\}.\) From property~\((b1),\)  we therefore have that the edges~\(f_{NW}(0,j)\) and~\(f_{NE}(0,j), m_1 = 1 \leq j \leq m_2-2, j \) odd, are consecutive edges of the cycle \(C_1.\) The edge \(f_{NE}(0,j)\) shares the endvertex \((0,j+1)\) with the edge \(f_{NW}(0,j+2)\) for all \(j\) and so~\(\Gamma_1\) is an subpath of the cycle~\(C_1.\)

It remains to see that every edge in \(\Delta_1 \setminus \{e_1,g_1\}\) has both endvertices in~\(R\setminus R_{left}.\) The edge \(e_1\) contains \((0,m_1-1)\) as an endvertex and the edge \(g_1\) contains \((0,m_{2}-1)\) as an endvertex. Also, every vertex \((0,j), 0 = m_1-1 \leq j \leq m_{2}-1\) either lies in the interior of the cycle~\(C_1\) or belongs to the path~\(\Gamma_1.\) Therefore if there exists an edge \(e \in Q_1\) intersecting the left edge \(R_{left}\) at some point \((0,y),\) then \(y \geq m_{2}\) or \(y \leq m_{1}-2.\)

If \(y \geq m_{2},\) then \(e\) is the edge of some square \(S_{0,y_1}\) with centre \((0,y_1),\) where \(y_1 \geq m_{2}.\) But this implies that \(y_1\) is odd and so \((0,y_1) \in {\cal C}\cap R_{left}.\) Thus the square~\(S_{0,y_1}\) is occupied and if it does not lie in the interior of the cycle~\(C_1,\) then we could merge \(S_{0,y_1}\) and \(C_1\) using Theorem~\ref{thm1} to form a bigger cycle containing the square \(S_{0,m_1}\) in its interior. This is a contradiction since square \(S_{0,m_1}\) is in the interior of the cycle~\(C_1\) and by construction, the cycle~\(C_1\) is the outermost boundary cycle containing the square \(S_{0,m_1}\) in its interior (see property \((c)\) of Lemma~\ref{outer}).

From the above paragraph, we therefore have that the square~\(S_{0,y_1}\) lies in the interior of the cycle~\(C_1\) and \(y_1 \geq m_2.\) But by definition, the integer \(m_2-2\) is the largest integer \(j\) such that \((0,j)\) lies in the interior of the cycle~\(C_1.\) Thus we have a contradiction and so every edge in \(\Delta_1 \setminus \{e_1,g_1\}\) belongs to \(R\setminus R_{left}.\) An analogous proof as above holds for \(y \leq m_1-2.\)

To prove the final statement of \((d1),\) suppose that some edge \(e\) in \(\Delta_1\) lies to the left of \(R_{left}.\) This necessarily means that \(e = e_1\) or \(e = g_1.\) We suppose \(e = g_1\) and arrive at a contradiction and an analogous proof holds for the other case. If \(e = g_1\)  lies to the left  of \(R_{left},\) then necessarily \(e = g_1 = f_{NW}(0,m_2).\) This is because \(g_1\) contains \((0,m_2-1)\) as an endvertex and the only edge lying to the left of \(R_{left}\) and containing \((0,m_2-1)\) as an endvertex is \(f_{NW}(0,m_2).\) But this means that the occupied square \(S_{0,m_2}\) containing \(f_{NW}(0,m_2)\) as an edge lies in the interior of the cycle \(C_1.\) In particular, the vertex \((0,m_2) \in {\cal C} \cap R_{left}\) lies
in the interior of the cycle \(C_1,\) a contradiction to the definition of the integer \(m_2\) (see paragraph following (\ref{m_def})). This proves~\((d1).\)

To prove property~\((d2),\) we proceed as follows. By the definition of the integers \(\{m_i\},\) we have that the cycles~\(C_{i}\) and~\(C_{i+1}\) intersect at the vertex~\((0,m_{i+1}-1) \in R_{left}.\) Using property~\((a3)\) of the previous subsection, we have that \(C_i\) and \(C_{i+1}\) intersect \emph{only} at \((0,m_{i+1}-1).\)

To prove the remaining part of \((d2),\) we proceed as follows. We prove the statement regarding the edge~\(f_{SW}(0,m_i).\) An analogous proof holds for the other statement regarding the edge~\(f_{SE}(0,m_{i+1}-2).\) Also we prove for \(i= 2\) and an analogous statement holds for all \(2 \leq i \leq h-1.\) The vertex \((0,m_{2}) \in {\cal C} \cap R_{left}\) lies in the interior of the cycle \(C_2\) and so we have from property~\((b1)\) that the edge~\(f_{NW}(0,m_{2})\) belongs to the cycle~\(C_{2}.\)

\begin{figure}[tbp]
\centering
\includegraphics[width=3.5in, trim= 50 350 50 150, clip=true]{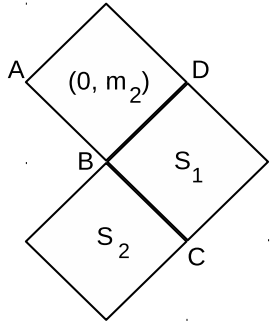}
\caption{The square \(S_{0,m_2}\) is represented by the centre \((0,m_2).\) The squares \(S_1\) and \(S_2\) represent the squares \(S_{1,m_2-1}\) and \(S_{0,m_2-2},\) respectively. }
\label{sm2_fig}
\end{figure}

In Figure~\ref{sm2_fig}, the square \(S_{0,m_2}\) is represented by its centre \((0,m_2)\) and the edge \(f_{NW}(0,m_2)\) is the line segment \(AB.\) If the edge \(f_{SW}(0,m_{2}) = BD\) does not belong to the cycle~\(C_{2},\) then it necessarily belongs to the interior of~\(C_2.\) Both the squares containing the edge \(f_{SW}(0,m_{2}) = BD\)  lie in the interior of the cycle \(C_2.\) The square~\(S_{1,m_{2}-1}\) with centre~\((1,m_2-1)\) (represented by \(S_1\) in Figure~\ref{sm2_fig}) also contains the edge~\(f_{SW}(0,m_2) = BD\) and therefore lies in the interior of \(C_2.\) 

The square~\(S_{0,m_2-2}\) with centre~\((0,m_2-2)\) (represented by \(S_2\) in Figure~\ref{sm2_fig}) shares the edge \(f_{NW}(1,m_2-1) = BC\) with the square~\(S_{1,m_2-1}.\) Since the square~\(S_{0,m_2-2}\) is occupied by definition, the square \(S_{0,m_2-2}\) also lies in the interior of the cycle~\(C_2.\)  We prove this by contradiction. The cycle \(C_2\) shares the edge \(f_{NW}(1,m_2-1) = BC\)  with the square \(S_{0,m_2-2}\) and so we can merge~\(C_2\) and~\(S_{0,m_2-2}\) using Theorem~\ref{thm1} and obtain a bigger cycle containing the square~\(S_{0,m_2}\) in its interior. This contradicts the fact that~\(C_2\) is the outermost boundary cycle containing the square \(S_{0,m_2}\) (see condition \((c)\) of Lemma~\ref{outer}).

From the above paragraph, we therefore have that the square~\(S_{0,m_2-2}\) lies in the interior of the cycle~\(C_2\) and by definition, the square~\(S_{0,m_2-2}\) also lies in the interior of the cycle \(C_1.\) Since the cycle~\(C_1\) and~\(C_2\) have mutually disjoint interiors (see property~\((a1)\)), we get a contradiction. This proves \((d2).\)

To prove \((d3),\) we consider the cycle graph \(H_{cyc}\) obtained as follows. Consider the vertex set \(\{1,\ldots,h\}\) and join \(i\) and \(j\) by an edge \(e(i,j)\) if the corresponding cycles \(C_i\) and \(C_j\) share a common vertex.

The graph~\(H_{cyc}\) constructed above is acyclic by property~(\ref{hcyc}) in the proof of Theorem~\ref{thm3}. Also, for every \(1 \leq i \leq h-1,\) the cycles \(C_i\) and \(C_{i+1}\) intersect and so the vertices \(i\) and \(i+1\) are by the edge \(e(i,i+1)\) in \(H_{cyc}.\) Therefore the path \((e(1,2),e(2,3),e(3,4),\ldots,e(h-1,h))\) is contained in the graph~\(H_{cyc}.\) If the cycle~\(C_2\) shares a vertex with the cycle~\(C_j\) for some \(j \neq 1,3\) then the vertices \(2\) and \(j\) would be also joined by an edge. This would mean that \(H_{cyc}\) contains a cycle, a contradiction. This proves~\((d3).\)~\(\qed\)

\subsubsection*{Outermost boundary with orientation}
We recall that \({\cal C}\) is the collection of vertices reachable by oriented open path starting from the left edge \(R_{left}\) of the rectangle \(R.\) We now introduce orientation only for the edges of the squares~\(S_z\) for vertices \(z \in {\cal C}\) and establish the properties needed for obtaining the oriented dual crossing. As before we denote the (dual) edges of the square~\(S_z\) as \(f_{NW}(z), f_{SW}(z), f_{NE}(z)\) and \(f_{SE}(z).\) We assign the orientations \(\nearrow, \nwarrow,\searrow\) and \(\swarrow,\) respectively, to the edges \(f_{NE}(z), f_{NW}(z), f_{SE}(z)\) and~\(f_{SW}(z)\) so that the edges in the square \(S_z\) form an oriented cycle.

For vertex \(z \in {\cal C},\) we recall that the corresponding square \(S_z\) is defined to be occupied. An edge belonging to \(S_z\) can also belong to another occupied square \(S_w\) and therefore have two possible orientations. However from property~\((a2),\) we have that every dual edge in the outermost boundary \(\cup_{i=1}^{h}C_i\) is a boundary edge adjacent to one occupied square and one vacant square. Therefore all dual edges in the outermost boundary have a unique orientation and we call them as \emph{oriented dual} edges.

Henceforth unless otherwise mentioned, we consider only oriented dual edges.


We have the following properties regarding the orientation of the cycles~\(C_i, 1 \leq i \leq h,\) of the outermost boundary. \\
\((f1)\) For each \(i, 1 \leq i \leq h,\) the subpath \((e_i,\Gamma_i,g_i) \subset C_i\) defined in \((n1)-(n2)\) is an oriented subpath of~\(C_i.\) Every edge in \(C_i \setminus (e_i,\Gamma_i,g_i)\) has both endvertices in \(R\setminus R_{left}.\)\\
\((f2)\) For each \(i, 1 \leq i \leq h,\) the cycle \(C_i\) is an oriented cycle; i.e., an oriented path with coincident endvertices.\\
\((f3)\) Suppose \(h \geq 2.\) For \(1 \leq i \leq h-1,\) the edge \(g_i = f_{SE}(0,m_{i+1}-2)\) has orientation \(\searrow.\) For \(2 \leq i \leq h,\) the edge \(e_i = f_{SW}(0,m_i)\) has orientation \(\swarrow.\)\\

\begin{figure}[tbp]
\centering
\includegraphics[width=3.5in, trim= 50 300 50 100, clip=true]{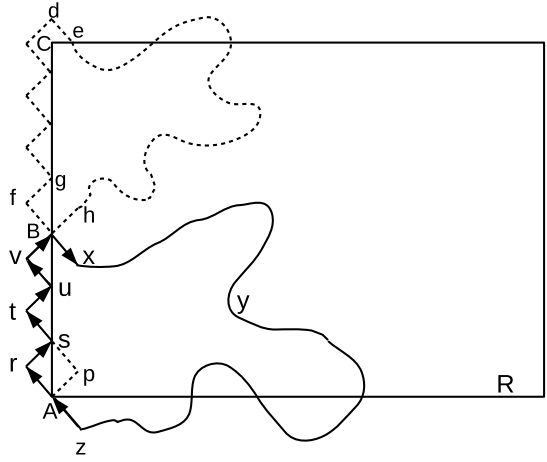}
\caption{The cycles \(C_1\) (thick line) and \(C_2\) (dotted line) are shown meeting at the point \(B\) on the left edge of the rectangle \(R.\) Also the oriented subpath \((e_1,\Gamma_1,g_1)\) of the cycle \(C_1\) is illustrated with endvertices~\(z\) and~\(x.\) }
\label{cyc_c1c2}
\end{figure}



The properties \((f1)\) and \((f3)\) are illustrated in Figure~\ref{cyc_c1c2}, where we have introduced orientation for the cycles~\(C_1\) and~\(C_2\) described in Figure~\ref{c1_c2_unor}. The oriented subpath \((e_1,\Gamma_1,g_1)\) is shown where the oriented edge \(e_1\) is denoted by the oriented segment~\(zA,\) the oriented path~\(\Gamma_1 = ArstuvB\) and the oriented edge~\(g_1 = Bx.\)


\emph{Proof of \((f1)-(f3)\)}: We prove \((f1)-(f2)\) and the first part of \((f3)\) for \(i = 1\) and an analogous proof holds for all \(i\) and the other cases.

The second statement of \((f1)\) is true by the final statement in property \((d1)\) and is used in the proof of \((f2).\) To see that the first statement of~\((f1)\) is true, we recall from property \((b1)\) that if \(z \in {\cal C} \cap R_{left},\) then the edges \(f_{NW}(z)\) and \(f_{NE}(z)\) of the square \(S_z\) are consecutive edges in some cycle \(C_j\) of the outermost boundary. Also, the square \(S_z\) is occupied and the square \(S_{z_1}\) sharing the edge \(f_{NW}(z)\) with \(S_z\) is vacant. Therefore \(f_{NW}(z)\) has the orientation \(\nwarrow.\) Similarly, the edge \(f_{NE}(z)\) has the orientation \(\nearrow.\) From property~\((d1),\) we have that \(\Gamma_1\) is a subpath of the cycle \(C_1.\) From the above, we obtain that~\(\Gamma_1\) is an \emph{oriented} subpath of the cycle~\(C_1.\)

The edge \(e_1\) defined in property \((d1),\) shares an endvertex with the edge \(f_{NW}(0,m_1) = f_{NW}(0,1)\) of the cycle \(C_1.\) Here we use \(m_1 = 1\) (see definition of \(m_i\) prior to the properties \((d1)-(d3)\)). There are therefore only two possibilities for~\(e_1.\) \((p1)\) Either \(e_1 = f_{SW}(0,1)\) (represented by the segment \(Ap\) in Figure~\ref{cyc_c1c2}) or \((p2)\) the edge \(e_1 = f_{NW}(1,0) = Az.\) For the case of \((p1),\) we argue as follows. From property \((a2),\) the edge \(e_1 = Ap\) is a boundary edge adjacent to one occupied square and one vacant square. The square~\(S_{0,1}\) (represented by \(AA_1A_2A_3\) in Figure~\ref{cyc_c1c2}) with centre \((0,1)\) (the midpoint of the segment~\(As\) in Figure~\ref{cyc_c1c2}) is occupied and is contained in the interior of the cycle \(C_1.\) Therefore the square~\(S_{1,0}\) with centre \((1,0)\) sharing the edge \(e_1 = f_{SW}(0,1) = Ap\) with the square \(S_{0,1}\) is vacant and lies in the exterior of~\(C_1.\) The edge \(e_1 = Ap\) therefore has unique orientation~\(\swarrow.\)

In the case of \((p2),\) the square \(S_{0,-1}\) with centre \((0,-1)\) sharing the edge~\(Az\) with \(S_{1,0}\) is vacant by definition since the vertex \((0,-1)\) lies in the exterior of the rectangle \(R.\) Therefore the square \(S_{1,0}\) with centre \((1,0)\) is occupied and so the edge \(e_1 = Az\) has unique orientation \(\nwarrow.\) In either case \((e_1,\Gamma_1)\) is an oriented subpath of cycle \(C_1.\) An analogous argument holds for the edge~\(g_1.\) This proves \((f1)\) and an analogous argument as above also proves the first statement of~\((f3).\) An analogous argument holds for the second statement in~\((f3).\)

To prove \((f2),\) we proceed by induction. Let \(C_1=  (q_1,\ldots,q_t)\) with \(\{q_i\}\) being the edges of \(C_1.\) From property \((f1),\) we assume that \(q_1 = f_{NW}(0,1),\)\\\(q_2 = f_{NE}(0,1),\)\( q_3 = f_{NW}(0,2),\ldots, q_{2m_2-3} = f_{NW}(0,m_{2}-2)\) and \(q_{2m_2-2} = f_{NE}(0,m_{2}-2).\) Also the last edge \(q_t = e_1\) and the edge \(q_{2m_2-1} = g_1\) so that \((e_1,\Gamma_1,g_1)\) is an oriented subpath of \(C_1.\) In Figure~\ref{cyc_c1c2}, the edge \(e_1 = Az\) and \(g_1 = Bx\) and the oriented subpath \((e_1,\Gamma_1,g_1)\) is shown with endvertices~\(z\) and~\(x.\)

\begin{figure*}
\centerline{\subfigure[]{\includegraphics[width=3.0in, trim= 100 350 100 150, clip=true]{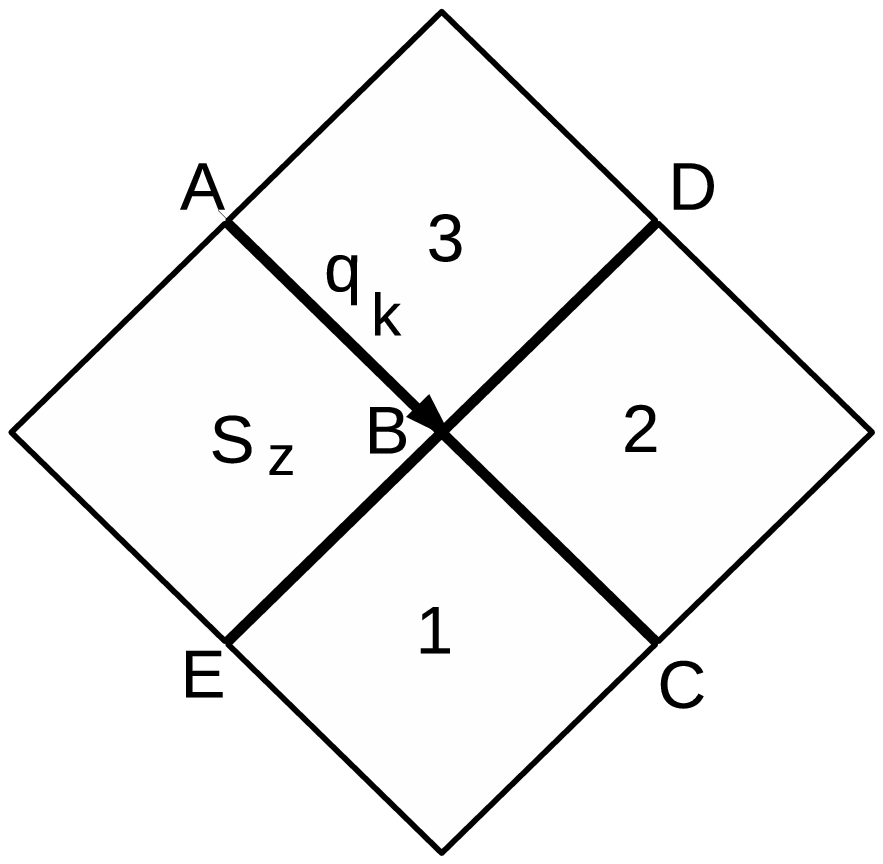}
\label{fig_first_case}}
\hfil
\subfigure[]{\includegraphics[width=3.0in, trim= 100 350 100 200,clip=true]{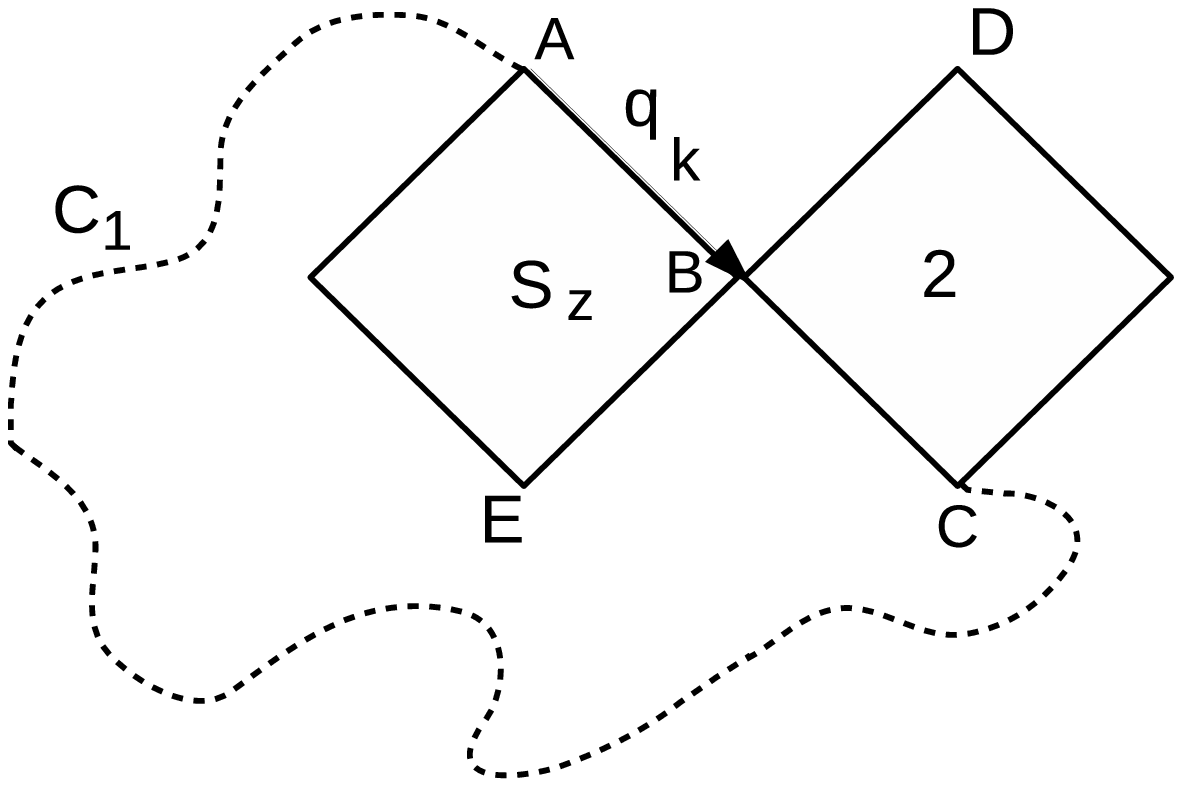}
\label{fig_second_case}}}
\caption{\((a)\) Squares neighbouring \(S_z.\) The label \(i\) corresponds to square \(S_{z_i}\) for \(i = 1,2,3.\) \((b)\) In case \((x3),\) exactly one of the squares~\(S_z\) or~\(S_{z_2}\) lies in the interior of~\(C_1.\)}
\label{sz_fig}
\end{figure*}

For the induction step, suppose \((q_t,q_1,q_2,\ldots,q_{k-1},q_k)\) is an oriented path for some \(2m_2-1 \leq k \leq t-2.\) We consider the case \(q_k = f_{SE}(z)\) for some vertex~\(z\) in the rectangle~\(R\) with orientation~\(\searrow.\)  An analogous proof holds for the other three types of oriented edges. We define the neighbouring squares of~\(S_z\) as follows. Let \(S_{z_1}\) and \(S_{z_3}\) be the squares sharing the edges \(f_{SE}(z)\)  and  \(f_{SW}(z),\) respectively with the square \(S_z.\) Let \(S_{z_2}\) be the square that shares a corner with \(S_{z}\) and share edges with \(S_{z_1}\) and \(S_{z_3}.\) This is illustrated in Figure~\ref{sz_fig}\((a)\) where the square \(S_z\) is shown along with the squares \(S_{z_i}, i = 1,2,3.\) The label \(i\) corresponds to the square \(S_{z_i},\) for \( i = 1,2,3.\) The oriented segment \(AB\) represents the edge \(q_k = f_{SE}(z).\)

We use the following property in the proof.
\begin{eqnarray}
&&\text{The square \(S_z\) is occupied and lies in the interior of the cycle \(C_1.\)}\nonumber\\
&&\;\;\;\;\;\;\text{The square~\(S_{z_3}\) is vacant and lies in the exterior of \(C_1.\)} \label{sz_prop}
\end{eqnarray}
\emph{Proof of (\ref{sz_prop})}: The edge \(q_k = f_{SE}(z) = f_{NW}(z_3)\) is common to both the squares~\(S_z\) and~\(S_{z_3}.\) From property~\((a2),\) we have that every edge in the cycle~\(C_1\) is a boundary edge adjacent to one occupied square contained in the interior of~\(C_1\) and one vacant square in the exterior of \(C_1.\) In particular, one of the squares \(S_z\) or \(S_{z_3}\) is vacant and the other is occupied. If \(S_{z_3}\) is occupied, then the orientation of \(q_k\) would be \(\nwarrow,\) a contradiction. Thus \(S_z\) is occupied and \(S_{z_3}\) is vacant and this proves~(\ref{sz_prop}).\(\qed\)

We now determine the endvertex common to the edges \(q_k\) and \(q_{k+1}.\) By induction assumption, we have that the edges \((q_{k-1},q_k)\) form a consistent pair and so the non arrow endvertex \(A\) of the edge~\(q_k\) (see Figure~\ref{sz_fig}) is coincident with the arrow endvertex of \(q_{k-1}.\) The vertex \(A\) has degree two in the cycle~\(C_1\) and so the edge \(q_{k+1}\) does not contain \(A\) as an endvertex. This means that the edge~\(q_k\) shares the arrow endvertex~\(B\) of the edge \(q_k\) and so the possible choices for the edge \(q_{k+1}\) are \(f_{SW}(z),f_{SE}(z_1)\) and \(f_{NE}(z_2),\) represented by the segments \(BE, BC\) and \(BD,\) respectively.

We consider three cases: \((x1)\) the squares \(S_{z_1}\) and \(S_{z_2}\)  are both vacant, \((x2)\) the square \(S_{z_1}\) is occupied and \(S_{z_2}\) is vacant and \((x3)\) the square \(S_{z_2}\) is occupied. In each case, we see that the edges \(q_k\) and \(q_{k+1}\) form a consistent pair.

Consider first the case~\((x1).\) The edge \(q_{k+1}\) is not \(f_{SE}(z_1) = BC\) since the latter edge is adjacent to two vacant squares~\(S_{z_1}\) and~\(S_{z_2}.\) We recall that every edge in the cycle \(C_1\) is a boundary edge belonging to an occupied square and a vacant square. Similarly, the edge~\(q_{k+1} \) is not \(f_{NE}(z_2) = BD\) since the latter is again adjacent to two vacant squares (see (\ref{sz_prop})). Thus \(q_{k+1} = f_{SW}(z) = BE\) with the corresponding orientation \(\swarrow.\)

The argument for case \((x2)\) is analogous as in the previous paragraph. The edge \(q_{k+1}\) cannot be \(f_{SW}(z) = BE\) since the latter belongs to two occupied squares~\(S_z\) and \(S_{z_1}\) (see (\ref{sz_prop})). Also \(q_{k+1}\) cannot be \(f_{NE}(z_2) = BD\) since the latter belongs to two vacant squares \(S_{z_2}\) and~\(S_{z_3}\) (again using \ref{sz_prop}). Therefore \(q_{k+1} = f_{SE}(z_1) = BC\) with the corresponding orientation~\(\searrow.\)

For case \((x3),\) we use the property \((a3)\) that the cycles \(\{C_i\}\) in the outermost boundary have mutually disjoint interiors and share at most one vertex in common. Moreover, the endvertex common to \(C_i\) and \(C_{i+1}\) belongs to the left edge \(R_{left}\) of \(R.\) The square~\(S_{z_2}\) is occupied  and shares a vertex \(B \notin R_{left}\) with the occupied square \(S_z\) contained in the interior of the cycle~\(C_1.\) The square~\(S_{z_2}\) therefore also belongs to the interior of the cycle \(C_1.\) Here \(B \notin R_{left}\) since the edge \(q_{k+1}\) belongs to the subpath \(C_1 \setminus (e_1,\Gamma_1,g_1)\) and therefore has both its endvertices in \(R \setminus R_{left}\) using property \((f1).\)

If the edge \(q_{k+1} = f_{NW}(z_2) = f_{SE}(z_1),\) then the cycle \(C_1\) encloses exactly one of the squares \(S_z\) or \(S_{z_2}\) in its interior. This is illustrated in Figure~\ref{sz_fig}\((b)\) where the segment \(AB\) is the edge \(q_k = f_{SE}(z)\) and the segment \(BC = f_{NW}(z_2).\) Therefore, \(q_{k+1} \neq f_{NW}(z_2)\) and so \(q_{k+1} = f_{NE}(z_2)\) with the orientation \(\nearrow.\) This completes the induction step for the case \(q_k = f_{SE}(z)\) for some vertex \(z \in {\cal C}\) with orientation \(\searrow.\) An analogous proof holds for the other three types of oriented edges.~\(\qed\)





\section{Proof of Theorem~\ref{thm8}}\label{pf8}
In this section, we prove the mutual exclusivity of open oriented left right and closed dual oriented top bottom dual crossings. Let \(LR_{or}(R)\) be the event that there exists an open oriented left right crossing and let \(TD^*_{or}(R)\) be the event that there is a closed dual oriented top bottom crossing of the rectangle~\(R.\) As in the unoriented case, the above two events cannot occur together. Suppose not and let \(\Pi_{L}\) be an open oriented left right crossing and let \(\Pi^*_D\) be a closed dual oriented top bottom crossing.

Every edge in \(\Pi_L\) lies in the interior of the rectangle \(R\) in the sense that if \(e \in \Pi_L\) then no endvertex of \(e\) lies in the exterior of~\(R.\) Also, since every edge in \(\Pi^*_D\) is closed, every edge in \(\Pi^*_D\) also lies in the interior of \(R\) (see the definition of closed dual edges two paragraphs prior to statement of Theorem~\ref{thm8}). Therefore the paths \(\Pi_L\) and \(\Pi^*_D\) intersect in the sense that there are  edges \(e \in \Pi_L\) and \(f \in \Pi^*_D\) that intersect. The edge \(e\) belonging to the original percolation model of Figure~\ref{or_mod} is open and the dual edge \(f\) intersecting \(e\) is closed, a contradiction to the definition of openness of the dual edges.

If the event \(LR_{or}(R)\) does not occur, we obtain a closed dual oriented top bottom crossing in a two step procedure described below.

\subsubsection*{Concatenating the paths \(\{\Delta_i\}_{1 \leq i \leq h}\) to form \(\Delta_{tot}\)}
The first step is to concatenate the paths \(\{\Delta_i\}\) defined in properties \((d1)-(d3).\) We recall the definition of the integers \(\{m_i\}\) defined in the previous subsection (see (\ref{m_def})). We also recall that \(m\) and \(n\) are the width and the height of the rectangle \(R.\) We have the following properties. \\
\((y1)\) Fix \(1 \leq i \leq h.\) The path \(\Delta_i\) has endvertices \(v_{i+1} = (0,m_{i+1}-1)\) and \((0,m_i-1).\) For \(i \neq j,\) the paths \(\Delta_i\) and \(\Delta_j\) are edge disjoint. At least one edge in the path~\(\Delta_h\) intersects the top edge \(R_{top}\) of the rectangle~\(R\) and at least one edge in the path~\(\Delta_{1}\) intersects the bottom edge \(R_{bottom}\) of~\(R.\)\\
\((y2)\) Let \(\Delta_{tot} = \cup_{i=1}^{h} \Delta_i.\) We have that
\(\Delta_{tot}\) is an oriented path with endvertices \((0,m_{h+1}-2)\) and \((0,m_1-1) = (0,0).\) If \(\Delta_{tot} = (h_1,\ldots,h_r),\) then the edge \(h_1 \in \Delta_{h}\) contains the endvertex \((0,m_{h+1}-2)\) and the edge~\(h_r \in \Delta_1\) contains the endvertex~\((0,0).\) Moreover there are indices \(k_1 < k_2\) such that~\(h_{k_1}\) intersects the top edge~\(R_{top}\) of the rectangle~\(R\) and~\(h_{k_2}\) intersects the bottom edge~\(R_{bottom}\) of~\(R.\)\\
\((y3)\) Let \(w = (x,y)\) be an endvertex of a (dual) edge \(e \in \Delta_{tot}.\) If \(w\) lies in the exterior of \(R,\) then the other endvertex of~\(e\) intersects the boundary of the rectangle~\(R.\) If \(R\) has no oriented left right crossing in the original percolation model (see Figure~\ref{or_mod}\((a)\)), then \(0 \leq x \leq m\) and \(-1 \leq y \leq n+1.\)\\

\emph{Proof of \((y1)-(y3)\)}: We prove \((y3), (y1)\) and \((y2)\) in that order. For~\((y3),\) we argue as follows. Let \(w = (x,y)\) be an endvertex of an edge \(e \in \Delta_i.\) If the vertex \(w\) lies in the interior of the rectangle~\(R,\) then \(0 \leq x \leq m\) and \( 0 \leq y \leq n.\)

If the vertex~\(w\) lies in the exterior of the rectangle~\(R,\) then the other endvertex \(w_1\) of the edge~\(e \in \Delta_i\) necessarily intersects the boundary of \(R.\) The endvertex \(w_1\) cannot lie in the interior of \(R\) because this means that \(w\) either belongs to the boundary or lies in the interior of \(R.\) If \(w_1\) does not belong to the boundary of \(R\) and lies in the exterior of \(R,\) then the square \(S_z\) with centre \(z\) and containing \(e\) as an edge lies in the exterior of \(R.\) This leads to a contradiction since all dual edges belong to squares \(S_z, z \in {\cal C}\) and no vertex in the oriented cluster \({\cal C}\) lies in the exterior of~\(R\) (see first subsection of Section~\ref{pf7}). Thus the endvertex \(w_1\) intersects the boundary of \(R\) and the bounds on the \(y-\)coordinate for the vertex \(w\) holds.

To obtain the bounds on the \(x-\)coordinate, we recall that \(e\) is an edge of the square \(S_z\) with centre \(z \in {\cal C}.\) Since the rectangle \(R\) has no left right crossing, the vertex \(z\) does not belong to the right edge \(R_{right}\) of~\(R.\) Thus both the endvertices of \(e\) lie to the left of \(R_{right}\) and so \(x \leq m.\) To see \(x \geq 0,\) we use property \((d1)\) that every edge in \(\Delta_i\) lies to the right of \(R_{left}.\) This proves the first statement of \((y3).\) The argument in the previous paragraph also proves the second statement of \((y3).\)

The first statement of \((y1)\) is true from property~\((d1)\) since the path \(\Delta_i\) has the same endvertices \((0,m_{i+1}-1)\) and \((0,m_i-1)\) as the path \(\Gamma_i = C_i \setminus \Delta_i.\)

For the second statement of \((y1),\) we use property \((a3)\) of the outermost boundary. The cycles \(C_i\) and \(C_j\) of the outermost boundary are edge disjoint for \(i \neq j.\) Since the path \(\Delta_i \subset C_i,\) we obtain the second statement of \((y1).\)

For the third statement of \((y1),\) we argue as follows. From the first statement of~\((y1),\) the path \(\Delta_1\) has \((0,m_1-1)\) as an endvertex. Since \(m_1 = 1\) (see (\ref{m_def})), the path~\(\Delta_1\) contains the origin as endvertex and so the second half of the third statement of~\((y1)\) is true.

To prove the first half, we again use the first statement of \((y1)\) that the path \(\Delta_h\) contains \((0,m_{h+2}-1)\) as an endvertex. Since \(m_{h+2}-1 \geq n\) (see~(\ref{m_h_1_def})), the endvertex \(v_{h+1}\) either lies above the top edge \(R_{top}\) or touches the top edge of the rectangle~\(R.\) Suppose endvertex \(v_{h+1}\) lies above \(R_{top}.\) Arguing as in the proof of~\((y3),\) we obtain that the other endvertex~\(w_{h+1}\) of the dual edge~\(g_h \in \Delta_h\) containing~\(v_{h+1}\) intersects the top edge~\(R_{top}.\) This proves the first half of the final statement of~\((y1)\) and therefore completes the proof of \((y1).\)

To prove \((y2),\) we have from property \((f2)\) that every \(\Delta_i, 1 \leq i \leq h\) is an oriented subpath of the cycle \(C_i.\) Therefore if \(h = 1,\) then \(\Delta_{tot} = \Delta_1\) is an oriented path. Suppose \(h = 2\) and consider the two subpaths \(\Delta_{1}\) and \(\Delta_{2}.\)

Using property \((d2),\) we have that \(\Delta_1\) and \(\Delta_{2}\) intersect only at the endvertex \((0,m_{2}-1).\) Thus \(\Delta_1 \cup \Delta_2\) is a path. To see that it \(\Delta_1 \cup \Delta_2\) an oriented path, we argue as follows. From property \((f3)\) we have that the last edge \(e_{2} \in \Delta_{2}\) is \(f_{SW}(0,m_{2})\) with orientation \(\swarrow.\) Also the first edge \(g_1 \in \Delta_1\) is \(f_{SE}(0,m_{2}-2)\) with orientation \(\searrow.\) From property \((d1)\) we have that both the edges \(e_{2}\) and \(g_1\) share the endvertex \((0,m_{2}-1).\) In particular, the arrow endvertex of \(e_{2}\) coincides with non arrow endvertex of \(g_{1}.\) Thus \((e_{2},g_1)\) is a consistent pair and we have that \(\Delta_1 \cup \Delta_{2}\) is an oriented path.

Suppose \(h \geq 3.\) Arguing as above, we have that \(\Delta_1 \cup \Delta_2\) is an oriented path. Consider now the union \(\Delta_1 \cup \Delta_2 \cup \Delta_3.\) Using property \((d2),\) we have that the path \(\Delta_3\) intersects the path \(\Delta_2\) at the point \((0,m_3-1).\) Arguing as in the previous paragraph, we have that \(\Delta_3 \cup \Delta_2\) is an oriented path. Using property \((d3),\) we also have that \(\Delta_3\) does not intersect the path \(\Delta_1.\) Thus \(\Delta_3 \cup \Delta_2\) intersects \(\Delta_1\) only at the vertex \((0,m_2-1)\) and so \(\Delta_3 \cup \Delta_2 \cup \Delta_1\) is also an oriented path. Continuing iteratively, we have that \(\Delta_{tot} = \cup_{i=1}^{h} \Delta_i\) is an oriented path with endvertices \((0,m_{h+1}-2)\) and \((0,m_1-1) = (0,0).\) This proves the first statement of~\((y2).\)

To prove the second statement of~\((y2)\)  we argue as follows. We have from property \((y1)\) that the path \(\Delta_h\) contains \((0,m_{h+1}-2)\) as an endvertex. Also no other path \(\Delta_i, 1 \leq i \leq h-1\) contains \((0,m_{h+1}-2)\) as an endvertex. Thus \(h_1 \in \Delta_h\) contains the endvertex \((0,m_{h+1}-2).\) Similarly the edge \(h_r \in \Delta_1\) contains the endvertex \((0,m_1-1)\) and since \(m_1 = 0\) (see (\ref{m_def})), this proves the second statement of \((y2).\) The final statement in~\((y2)\) follows from property \((y1).\)~\(\qed\)

\subsubsection*{Extracting the dual crossing from \(\Delta_{tot}\)}
Let \(\Delta_{tot} = (h_1,\ldots,h_r).\) Using property \((y2)\) above, let \(j_1\) be the largest index \(j < k_2\) such that the edge \(h_j\) intersects the top edge of \(R.\) Similarly, again using property \((y2),\) let \(j_2\) be the smallest index \(j > j_1\) such that \(h_j\) intersects the bottom edge of~\(R.\)
We have the following properties.\\
\((z1)\) The edge \(h_{j_1}\) touches the top edge of \(R\) and has orientation \(\searrow \) or \(\swarrow.\)\\
\((z2)\) The edge \(h_{j_2}\) touches the bottom edge of \(R\) and has orientations~\(\searrow\) or~\(\swarrow.\)\\
\((z3)\) If \(R\) has no oriented crossing in the original oriented percolation model of Figure~\ref{or_mod}\((a),\) then every edge in \(\Delta_{cr} = (h_{j_1},h_{j_1+1},\ldots,h_{j_2})\) lies in the interior of \(R.\)\\
Thus the path~\(\Delta_{cr}\) is the desired dual oriented top bottom crossing.


\emph{Proof of \((z1)-(z3)\)}: We prove \((z3)\) first. Suppose some edge \(h_k \in \Delta_{cr}, k \geq j_1\) lies in the exterior of the rectangle~\(R.\) From property~\((y3)\) above, we have that at least one endvertex of \(h_k\) intersects the boundary of \(R.\) Suppose \(h_k\) intersects the top edge of~\(R.\) We arrive at a contradiction as follows. Let \(v_{k}\) and \(v_{k+1}\) be the non arrow and arrow endvertices of the edge~\(h_{k},\) respectively. Suppose that one of \(v_{k}\) or \(v_{k+1}\) lies in the exterior of the rectangle~\(R.\) If the non arrow endvertex \(v_{k}\) lies in the exterior of \(R,\) the edge~\(h_{k+1}\) sharing the non arrow endvertex \(v_{k}\) also intersects the top edge of~\(R.\) This contradicts the definition of the index~\(j_1.\)

If the arrow endvertex \(v_{k+1}\) of the edge \(h_k\) lies in the exterior of \(R,\) then the edge \(h_{k+1}\) also lies in the exterior of \(R\) and the arrow endvertex of \(h_{k+1}\) intersects the top edge of \(R.\) This again contradicts the definition of \(j_1\) and so the edge \(h_k\) does not touch the top edge of \(R.\)


An analogous proof holds if \(h_k\) intersects the bottom edge of \(R\) and lies in the exterior of \(R.\) If the edge \(h_k\) intersects the left edge \(R_{left}\) of \(R\) and lies in the exterior of \(R,\) then \(h_k\) lies to the left of \(R_{left}.\) But since \(h_k \in \cup_{i=1}^{h}\Delta_i,\) this contradicts the final statement of property \((d1).\) Finally suppose \(h_k\) lies to the right of \(R_{right}\) and let \(x_k\) and \(x_{k+1}\) be the \(x-\)coordinates of the vertices \(v_k\) and \(v_{k+1},\) respectively. We must have either \(x_k \geq m+1\) or \(x_{k+1} \geq m+1,\) a contradiction to property \((y3).\) Thus every edge in \(\Delta_{cr}\) lies in the interior of \(R\) and this proves .

Using an argument analogous to the first two paragraphs of proof of \((z1)\) we obtain that \(h_{j_1}\) has orientation \(\searrow\) or \(\swarrow.\) This proves \((z1)\) and an analogous proof holds for \(h_{j_2}\) in the property~\((z2).\) \(\qed\)

\emph{Proof of Theorem~\ref{thm8}}: It only remains to see that every edge in \(\Delta_{cr}\) is closed and we prove this as follows. Let \(e\) be an oriented dual edge and let \(S_{a}\) and \(S_b\) be the squares with centres \(a\) and \(b,\) respectively, containing the edge \(e.\) Here \(a,b\) are vertices in the rectangle \(R.\) The edge \(e\) belongs to some cycle \(C_i\) of the outermost boundary and so from property \((a2),\) we have that one of the squares in \(\{S_a,S_b\}\) is occupied and the other is vacant. We recall that \(S_z, z \in R\) is occupied if the vertex \(z\) belongs to the oriented cluster \({\cal C}\) defined in the second paragraph of Section~\ref{pf7}. Thus one of the vertices of \(\{a,b\}\) belongs to \({\cal C}\) and the other does not belong to \({\cal C}.\) Therefore the edge \(f\) joining \(a\) and \(b\) is closed and so the edge \(e\) intersecting the edge \(f\) is also closed (see two paragraphs prior to statement of Theorem~\ref{thm8}).~\(\qed\)


\renewcommand{\theequation}{\thesection.\arabic{equation}}

\subsection*{Acknowledgement}
I thank Professors Rahul Roy and Federico Camia for crucial comments and for my fellowships. I also thank NISER for my fellowship.

\bibliographystyle{plain}

\end{document}